\documentclass{amsart}

\RequirePackage{color}

\def\smallskip{\vskip\smallskipamount}
\def\medskip{\vskip\medskipamount}
\def\bigskip{\vskip\bigskipamount}

\usepackage{amsmath,amsthm,bm,mathrsfs,amscd,tikz}
\usepackage{ytableau}
\usepackage{comment}
\usepackage{mathdots}
\usepackage[all]{xy}
\usepackage{graphicx}
\usepackage[colorinlistoftodos]{todonotes}
\usepackage{enumerate}
\usepackage{feynmp-auto}
\usepackage[english]{babel}
\usepackage[utf8]{inputenc}
\usepackage{float}
\usepackage{tikz}
\usepackage{tikz-cd}
\usepackage{amssymb}
\usepackage{mathtools}
\usetikzlibrary{decorations.pathmorphing}

\usepackage[utf8]{inputenc}

\newtheoremstyle{thmstyle}{}{}{\itshape}{}{\bfseries}{ }{5pt}{}
\newtheoremstyle{exstyle}{}{}{}{}{\bfseries}{ }{5pt}{}
\newtheoremstyle{defstyle}{}{}{}{}{\bfseries}{ }{5pt}{}
\newtheoremstyle{remstyle}{}{}{}{}{\bfseries}{ }{5pt}{}

\theoremstyle{thmstyle}
\newtheorem{thm}{Theorem}[section]
\newtheorem{theorem}[thm]{Theorem}

\newtheorem{lemma}[thm]{Lemma}

\newtheorem{proposition}[thm]{Proposition}

\newtheorem{corollary}[thm]{Corollary}

\newtheorem{conjecture}{Conjecture}[section]
\newtheorem{question}{Question}[section]

\theoremstyle{exstyle}

\newtheorem{example}[thm]{Example}

\theoremstyle{defstyle}

\newtheorem{definition}[thm]{Definition}

\newtheorem{def-prop}[thm]{Definition-Proposition}
\newtheorem{def-lem}[thm]{Definition-Lemma}
\newtheorem{rem-convention}[thm]{Remark-Convention}
\newtheorem{def-note}[thm]{Definition-Notation}

\theoremstyle{remstyle}

\newtheorem{remark}[thm]{Remark}
\newtheorem{construction}[thm]{Construction}
\newtheorem{Setting}[thm]{Setting}

\theoremstyle{remstyle}

\newcommand{\Hom}{\operatorname{Hom}}
\newcommand{\Fac}{\operatorname{Fac}}

\DeclareMathOperator*{\rad}{rad}

\DeclareMathOperator*{\modu}{mod}
\DeclareMathOperator*{\Modu}{Mod}

\newcommand{\Str}{\operatorname{Str}}

\DeclareMathOperator*{\Brick}{Brick}
\DeclareMathOperator*{\ind}{ind}
\DeclareMathOperator*{\Ind}{Ind}
\DeclareMathOperator*{\tors}{tors}

\DeclareMathOperator*{\ftors}{f-tors}

\DeclareMathOperator*{\End}{End}

\DeclareMathOperator*{\brick}{brick}

\DeclareMathOperator*{\rigid}{rigid}
\DeclareMathOperator*{\irigid}{\mathtt{i}rigid}
\DeclareMathOperator*{\ttilt}{-tilt}
\DeclareMathOperator*{\tilt}{tilt}
\DeclareMathOperator*{\trig}{-rigid}

\DeclareMathOperator*{\node}{node}

\DeclareMathOperator*{\pd}{pd}

\DeclareMathOperator*{\Mri}{Mri}

\DeclareMathOperator*{\nD}{nD}

\DeclareMathOperator*{\sB}{sB}

\DeclareMathOperator*{\GL}{GL}

\DeclareMathOperator*{\rep}{rep}

\DeclareMathOperator*{\Irr}{Irr}

\newcommand{\cupdot}{\mathbin{\mathaccent\cdot\cup}}

\makeatletter
\newcommand{\doublewidetilde}[1]{{%
  \mathpalette\double@widetilde{#1}%
}}
\newcommand{\double@widetilde}[2]{%
  \sbox\z@{$\m@th#1\widetilde{#2}$}%
  \ht\z@=.9\ht\z@
  \widetilde{\box\z@}%
}
\makeatother

\begin{document}

\title{Biserial algebras and generic bricks}
\author[Kaveh Mousavand, Charles Paquette]{Kaveh Mousavand, Charles Paquette} 
\address{Representation Theory and Algebraic Combinatorics Unit, Okinawa Institute of Science and Technology (OIST), Japan}
\email{mousavand.kaveh@gmail.com}
\address{Department of Mathematics and Computer Science, Royal Military College of Canada, Kingston ON, Canada}
\email{charles.paquette.math@gmail.com}
\thanks{The second-named author was supported by the National Sciences and Engineering Research Council of Canada (RGPIN-2018-04513), and by the Canadian Defence Academy Research Programme.}


\maketitle

\begin{abstract}
We use generic bricks in the study of arbitrary biserial algebras. For a biserial algebra $\Lambda$ of rank $n$ over an algebraically closed field $k$, we show that $\Lambda$ is brick-infinite if and only if it admits a generic brick, which we further prove to be the case if and only if there exists an infinite family of bricks of length $d$, for some $2\leq d\leq 2n$. 
Consequently, we obtain an algebro-geometric realization of $\tau$-tilting finiteness of biserial algebras: $\Lambda$ is $\tau$-tilting finite if and only if, for each dimension vector $\underline{d}$, there are only finitely many orbits of bricks  in the representation variety $\modu (\Lambda, \underline{d})$.
Our results rely on our full classification of minimal brick-infinite biserial algebras in terms of quivers and relations, seen as the modern analogue of the classification of minimal representation-infinite (special) biserial algebras, given by Ringel. We show that if $\Lambda$ is a minimal brick-infinite biserial algebra, then $\Lambda$ is gentle and admits exactly one generic brick. In this case, we describe the spectrum of $\Lambda$ and prove that it is similar to that of a tame hereditary algebra. In other words, $\Brick(\Lambda)$ is the disjoint union of a unique generic brick with a countable infinite set of bricks of finite lengths, and a family of bricks of length $d$ parameterized by $k$.
Our work strengthens and generalizes some earlier results on $\tau$-tilting finiteness of gentle algebras and special biserial algebras, respectively treated by Plamondon and Schroll-Treffinger-Valdivieso.
\end{abstract}

\tableofcontents
\section{Introduction}\label{Section:Introduction}
Throughout, all algebras are assumed to be finite dimensional associative unital over algebraically closed field. By $\Lambda$ we denote any such algebra and, with no loss of generality, we can further assume every algebra is basic and connected.
Hence, $\Lambda$ has a presentation of the form $kQ/I$, for a unique finite and connected quiver $Q$ and an admissible ideal $I$ in the path algebra $kQ$. Consequently, $\Lambda$-modules can be identified with representations of the bound quiver $(Q,I)$. Unless specified otherwise, we always work with left $\Lambda$-modules and consider them up to isomorphism. In particular, $\Modu \Lambda$ denotes the category of all left $\Lambda$-modules, whereas $\modu \Lambda$ denotes the category of all finitely generated left $\Lambda$-modules. 
Moreover, let $\Ind(\Lambda)$ and $\ind(\Lambda)$ respectively denote the collections of (isomorphism classes) of indecomposable modules in $\Modu \Lambda$ and $\modu \Lambda$.
The standard notations and terminology which are not explicitly defined here can be found in \cite{ASS}, or else they will be introduced throughout the text.

\subsection{Motivations and background} \label{Subsection:Motivations and background}
Recall that $M$ in $\Modu \Lambda$ is a \emph{brick} if $\End_{\Lambda}(M)$ is a division algebra. Then, $\Lambda$ is called \emph{brick-finite} if it admits only finitely many bricks (up to isomorphism). Each brick is evidently indecomposable and by $\Brick(\Lambda)$ and $\brick(\Lambda)$ we denote the subsets of all bricks, respectively in $\Ind(\Lambda)$ and $\ind(\Lambda)$.
Although each representation-finite (rep-finite, for short) algebra is brick-finite, the converse is not true in general (e.g. any representation-infinite local algebra admits a unique brick). 
More precisely, the notion of brick-finiteness is of interest only if $\Lambda$ is a rep-infinite tame algebra, or else when $\Lambda$ is wild but not strictly wild (a standard argument yields that any strictly wild algebra is brick-infinite).

Bricks and their properties play pivotal roles in different areas and can be studied from various perspectives (see Section \ref{Section:Preliminaries}).
As we do here, the notion of brick-finiteness can be viewed as a conceptual counterpart of representation-finiteness. To better highlight this perspective, let us present two analogous characterizations that also motivate our work. 
Thanks to some classical and recent results, $\Lambda$ is known to be rep-finite if and only if $\Ind(\Lambda)= \ind(\Lambda)$, whereas it is brick-finite if and only if $\Brick(\Lambda)= \brick(\Lambda)$ (for details, see \cite{Se}). Furthermore, through the lens of approximation theory, $\Lambda$ is shown to be rep-finite if and only if every full subcategory of $\modu \Lambda$ is functorially finite, whereas it is brick-finite if and only if every torsion class in $\modu \Lambda$ is functorially finite (for details, see \cite{DIJ}). 

Before we recall a powerful tool in the study of rep-infinite algebras, observe that each $\Lambda$-module $M$ can be viewed as a (right) module over $\End_{\Lambda}(M)$. Then, the \emph{endolength} of $M$ is the length of $M$ when considered as an $\End_{\Lambda}(M)$-module. In particular, a $\Lambda$-module $G$ in $\Ind(\Lambda)\setminus \ind(\Lambda)$ is called \emph{generic} if it is of finite endolength. 
Generic modules are known to play a significant role in representation theory of algebras. For instance, any generic $\Lambda$-module of endolength $d$ gives rise to an infinite family of (non-isomorphic) modules of length $d$ in $\ind(\Lambda)$ (see \cite{CB1} and the references therein). In fact, Crawley-Boevey \cite{CB1} has given an elegant realization of the Tame/Wild dichotomy theorem of Drozd \cite{Dr} in terms of generic modules and their endolength (see Theorem \ref{Thm Crawley-Boevey tame/domestic alg.}).
Based on his characterization, one can further refine tame algebras and say $\Lambda$ is \emph{$m$-domestic} if it admits exactly $m$ generic modules (up to isomorphism). In general, $\Lambda$ is \emph{domestic} if $\Lambda$ is $m$-domestic for some $m \in \mathbb{Z}_{\geq 0}$.

For an algebra $\Lambda$, it is known that $\Lambda$ is rep-infinite if and only if it admits a generic module (\cite{CB1}). Hence, it is natural to ask for an analogous characterization of brick-infinite algebras. 
First, observe that if $\Lambda$ is brick-finite, no generic $\Lambda$-module is a brick (see Theorem \ref{Thm:Sentieri}).
To describe our proposal for a new treatment of brick-infinite case, we introduce some new terminology.
In particular, we say $G$ is a \emph{generic brick} of $\Lambda$ if $G$ is a generic module and it belongs to $\Brick(\Lambda)$.
Furthermore, $\Lambda$ is called \emph{brick-continuous} if for a positive integer $d$, there exists an infinite family in $\brick(\Lambda)$ consisting of bricks of length $d$. The latter notion is motivated by the algebro-geometric properties of bricks (see Section \ref{Section:Some applications and problems}). In particular, $\Lambda$ is said to be \emph{brick-discrete} if it is not brick-continuous. We remark that brick-discrete algebras are studied in \cite{CKW}, where the authors call them ``Schur representation-finite" algebras and view them as a generalization of rep-finite algebras. In \cite{Mo2}, it is conjectured that brick-discrete algebras are the same as brick-finite algebras. 

\subsection{Problem and results}\label{Subsection:Problem and results}
In the rest of this section, we focus on our main problem and present our results. 
In particular, our conjecture below is primarily inspired by our earlier results in \cite{Mo1, Mo2} on the behavior of bricks, as well as our work on minimal $\tau$-tilting infinite algebras in \cite{MP} (see also Conjecture \ref{My Old Conjecture}).

\begin{conjecture}[Conjecture \ref{Conjecture for arbitrary algs}]\label{Introduction: Conjecture}
For an algebra $\Lambda$, the following are equivalent:
\begin{enumerate}
    \item $\Lambda$ is brick-infinite;
    \item $\Lambda$ is brick-continuous;
    \item $\Lambda$ admits a generic brick.
\end{enumerate}
\end{conjecture}

To verify the above conjecture in full generality, one only needs to treat those algebras which are brick-infinite and minimal with respect to this property. In particular, we say $\Lambda$ is a \emph{minimal brick-infinite} algebra (min-brick-infinite, for short) if $\Lambda$ is brick-infinite but all proper quotient algebras of $\Lambda$ are brick-finite. Thus, a concrete classification of such algebras will be helpful in the study of our conjecture. 
As discussed in Subsection \ref{Subsection:Minimal brick-infinite algebras}, min-brick-infinite algebras are novel counterparts of minimal representation-infinite algebras (min-rep-infinite, for short). This classical family is extensively studied due to their role in several fundamental problems, particularly in the celebrated Brauer-Thrall Conjectures (see \cite{Bo1} and references therein).
Note that min-brick-infinite algebras also enjoy some important properties that could be helpful in the study of Conjecture 1.1 (see Theorem \ref{min-tau-inf properties}).

It is known that biserial algebras form an important family of tame algebras, among which string, gentle and special biserial algebras have appeared in many areas of research (for definitions and details, see Subsection \ref{Subsection:Biserial algebras}). In particular, if $\Lambda$ is a special biserial algebra, one can combinatorially describe all indecomposable modules in $\modu \Lambda$, as well as their Auslander-Reiten translate and morphisms between them (for details, see \cite{BR} and \cite{WW}). However, for arbitrary biserial algebras, there is no full classification of indecomposable modules and their representation theory is more complicated than that of special biserial algebras.

In 2013, Ringel \cite{Ri1} gave an explicit classification of those min-rep-infinite algebras which are special biserial. Thanks to the more recent results of Bongartz \cite{Bo1}, one can show that Ringel's classification is in fact the full list of min-rep-infinite biserial algebras. In particular, any such biserial algebra is a string algebra which is either $1$-domestic or else non-domestic (for more details, see \cite{Ri1} and Subsection \ref{Subsection:Minimal brick-infinite algebras}). 
Here, we obtain an analogue of Ringel's classification and give a full list of min-brick-infinite biserial algebras in terms of their quivers and relations. 

Due to some technical observations that are further explained in Section \ref{Section:Some applications and problems}, in this paper, unless stated otherwise, we restrict to tame algebras and verify the above conjecture for all biserial algebras.
Before stating our main classification result, let us fix some new terminology. In particular, for an arbitrary algebra $\Lambda$, we say $\Lambda$ is \emph{$m$-generic-brick-domestic} if it admits exactly $m$ (isomorphism classes) of generic bricks. More generally, $\Lambda$ is \emph{generic-brick-domestic} if it admits only finitely many generic bricks.
For the definition and concrete description of \emph{generalized barbell algebras}, we refer to Subsection \ref{Subsection:Minimal brick-infinite algebras} and Figure \ref{fig:Generalized barbell quiver}. In particular, the following theorem follows from our results in Section \ref{Section:Minimal brick-infinite biserial algebras}.

\begin{theorem}\label{Thm: Introduction: min-brick-inf biserial}
Let $\Lambda$ be a minimal brick-infinite biserial algebra. Then, $\Lambda$ is either a hereditary algebra of type $\widetilde{A}_n$ or $\Lambda$ is a generalized barbell algebra. In particular, $\Lambda$ is a $1$-generic-brick-domestic gentle algebra.
\end{theorem}

If $Q$ is an affine Dynkin quiver, the algebra $kQ$ is known to be minimal brick-infinite and tame. 
A classical result of Ringel \cite{Ri2} shows that every such path algebra admits a unique generic brick. Hence, from this point of view, the preceding theorem extends the aforementioned result of Ringel and treats some non-hereditary tame minimal brick-infinite algebras (see Theorem \ref{Thm:Generarlization of Ringel's 1-generic-brick-domestic}).  

We remark that the generalized barbell algebras are never domestic. However, by the above theorem, they are always $1$-generic-brick-domestic. As an interesting consequence of our classification, the following corollary is shown in Section \ref{Section:Some applications and problems}.

\begin{corollary}\label{Cor: m-domestic and n-gen-brick-domestic}
For any $m$ and $n$ in $\mathbb{Z}_{\geq 0}$, there exists a string algebra $A$ where
\begin{enumerate}
    \item $A$ is $m$-domestic and $n$-generic-brick-domestic, with $m\geq n$.
    
    \item $A$ is $n$-generic-brick-domestic but of non-polynomial growth.
\end{enumerate}
Furthermore, there exist string algebras which are not generic-brick-domestic.
\end{corollary}

To prove the above corollary, in Section \ref{Subsection: Domestic vs. generic-brick-domestic} we present several examples and as the result give an explicit algorithm to construct an $n$-generic-brick-domestic algebra for each $n \in \mathbb{Z}_{\geq 0}$. Moreover, in Example \ref{Example: gentle algebra not generic-brick-domestic}, we give a gentle algebra which is not generic-brick-domestic. 
The above corollary further highlights the fundamental differences between generic modules and generic bricks, as well as the domestic and generic-brick-domestic algebras (for example, see Questions \ref{Question:generic-brick-domesticness of tame min-brick-infinite algs}).

As an important consequence of our classification result, we obtain a useful characterization of brick-finiteness of biserial algebras.
Note that the rank of the Grothendieck group of $\Lambda=kQ/I$ is the number of vertices in $Q$, denoted by $|Q_0|$.

\begin{theorem}\label{Thm:Introduction:brick-infinite biserial algs}
For a biserial algebra $\Lambda=kQ/I$, the following are equivalent:
\begin{enumerate}
\item $\Lambda$ is brick-infinite;
\item For some $ 2\leq d \leq 2|Q_0|$, there is an infinite family $\{X_{\lambda}\}_{\lambda \in k^*}$ in $\brick(\Lambda)$ with $\dim_{k}(X_{\lambda})= d$;
\item $\Lambda$ admits a generic brick whose endolength is at most $2|Q_0|$;
\end{enumerate}
\end{theorem}

The above theorem asserts a stronger version of Conjecture \ref{Introduction: Conjecture} for the family of biserial algebras.
Moreover, the numerical condition given in part $(2)$ is similar to that of Bongartz's for the length of $1$-parameter families of indecomposable modules over rep-inf algebras (for details, see \cite{Bo2}). This opens some new directions in the study of distribution of bricks (for example, see Question \ref{Question: Numerical test for brick-finiteness}).

As mentioned earlier, any systematic study of bricks and their properties provide new insights into several other domains of research. We end this section by the following corollary which highlight these connections and postpone further applications of our results to our future work. 
All the undefined terminology and notations used in the next assertion appear in Sections \ref{Section:Preliminaries} or \ref{Section:Minimal brick-infinite biserial algebras}. Moreover, proof of the following corollary follows from Theorem \ref{Thm:tau-tilting finiteness equivalences} and Corollary \ref{Corollary: all equivalences for brick-inf biserial algs}.

\begin{corollary}\label{Corollary-Introduction: all equivalences for brick-inf biserial algs}
Let $\Lambda$ be a biserial algebra. The following are equivalent:

\begin{enumerate}
    \item $\Lambda$ is $\tau$-tilting infinite;
    \item There is a band component $\mathcal{Z}$ in $\Irr(\Lambda)$ which contains a rational curve $\mathcal{C}$ of non-isomorphic bricks $\{M_{\lambda}\}$ such that $\mathcal{Z}=\overline{\bigcup_{\lambda\in \mathcal{C}}\mathcal{O}_{M_{\lambda}}}$;
    \item For some $\theta \in K_0(A)$, there exist infinitely many non-isomorphic $A$-modules which are $\theta$-stable.
\end{enumerate}
\end{corollary}

\section{Preliminaries}\label{Section:Preliminaries}
In this section we mainly collect some essential tools used in our main arguments. For the well-known results on (special) biserial algebras, tilting and $\tau$-tilting theory, as well as for the rudiments of representation varieties, we only provide references.

\subsection{Notations and conventions}
By a \emph{quiver} we always mean a finite directed graph, formally given by a quadruple $Q=(Q_0,Q_1,s,e)$, with the vertex set $Q_0$ and arrow set $Q_1$, and the functions $s,e: Q_1 \rightarrow Q_0$ respectively send each arrow $\alpha$ to its start. For $\alpha$ and $\beta$ in $Q_1$, by $\beta \alpha$ we denote the path of length two which starts at $s(\alpha)$ and ends at $e(\beta)$. 
Let  $Q_1^{-1}:=\{\gamma^{-1} \,| \,\gamma \in Q_1\}$ be the set of formal inverses of arrows of $Q$. That is, $s(\gamma^{-1})=e(\gamma)$ and $e(\gamma^{-1})=s(\gamma)$.

Following our assumptions in Section \ref{Section:Introduction}, every algebra $\Lambda$ is an admissible quotient $kQ/I$ of a path algebra $kQ$ for some quiver $Q$, up to Morita equivalence. In this case, the pair $(Q,I)$ is called \emph{a bound quiver}. All quotients of path algebras will be assumed to be admissible quotients. Moreover, modules over $\Lambda$ can be seen as representations over the bound quiver $(Q,I)$.
Provided we begin from a bound quiver, this dictionary is still available and $k$ can be an arbitrary field. In this case, $\Lambda$-modules are representations of the corresponding bound quiver. For $M$ in $\modu \Lambda$, let $|M|$ denote the number of non-isomorphic indecomposable modules that appear in the Krull-Schmidt decomposition of $M$. In particular, for $\Lambda=kQ/I$, we have $|\Lambda|=|Q_0|$, which is the same as the rank of $K_0(\Lambda)$, where $K_0(\Lambda)$ denotes the Grothendieck group of $\modu \Lambda$.
In particular, this rank is the number of (isomorphism classes of) simple modules in $\modu \Lambda$. 

For $\Lambda=kQ/I$, unless specified otherwise, we consider a minimal set of uniform relations that generate the admissible idea $I$. That is, each generator of $I$ is a linear combination of the form $R=\sum_{i=1}^{t} \lambda_i p_i$, where $t\in \mathbb{Z}_{>0}$ and $\lambda_i \in k\setminus \{0\}$, and all $p_i$ are paths of length strictly larger than one in $Q$ starting at the same vertex $x$ and ending at the same vertex $y$.
For the most part, we work with monomial and binomial relations, which are respectively when $t=1$ and $t=2$. In particular, the monomial relations of length $2$, known as quadratic monomial relations, play a crucial role in the study of (special) biserial algebras. A vertex $v$ in $Q$ is a \emph{node} if it is neither a sink nor a source, and for any arrow $\alpha$ incoming to $v$ and each arrow $\beta$ outgoing from $v$, we have $\beta\alpha \in I$.

In this paper, all subcategories are assumed to be full and closed under isomorphism classes, direct sum and summands.
Moreover, for a given collection $\mathfrak{X}$, we say a property holds for \emph{almost all} elements in $\mathfrak{X}$ if it is true for all but at most finitely many elements of $\mathfrak{X}$.

\subsection{Biserial algebras}\label{Subsection:Biserial algebras}
An algebra $\Lambda$ is said to be \emph{biserial} if for each left and right indecomposable projective $\Lambda$-module $P$, we have $\rad(P)=X+Y$, where $X$ and $Y$ are uniserial modules and $X\cap Y$ is either zero or a simple module. Biserial algebras were formally introduced by Fuller \cite{Fu}, as a generalization of uniserial algebras, and Crawley-Boevey \cite{CB2} showed that they are always tame.

Special biserial algebras form a well-known subfamily of biserial algebras and thanks to their rich combinatorics, their representation theory is well-studied. We recall that an algebra $\Lambda$ is \emph{special biserial} if it is Morita equivalent to an algebra $kQ/I$ such that the bound quiver $(Q,I)$ satisfies the following conditions:
\begin{enumerate}[(B1)]
\item At every vertex $x$ in $Q_0$,  there are at most two incoming and at most two outgoing arrows.
\item For each arrow $\alpha$ in $Q_1$, there is at most one arrow $\beta$ such that $\beta \alpha \notin I$ and at most one arrow $\gamma$ such that $\alpha \gamma \notin I$.
\end{enumerate}
 
A special biserial algebra $\Lambda=kQ/I$ with $(Q,I)$ as above is called a \emph{string algebra} if $I$ in $kQ$ can be generated by monomial relations. 
Over string algebras, all indecomposable modules and morphisms between them are understood (see \cite{BR} and \cite{WW}).

An important subfamily of string algebras consists of gentle algebras. Recall that $\Lambda=kQ/I$ is \emph{gentle} if it is a string algebra and $I$ can be generated by a set of quadratic monomial relations such that $(Q,I)$ satisfies the following condition:

\begin{enumerate}
\item[(G)] For each arrow $\alpha \in Q_1$, there is at most one arrow $\beta$ and at most one arrow $\gamma$ such that $0\ne\alpha\beta \in I$ and $0\ne\gamma\alpha \in I$.
\end{enumerate}

Observe that if $\Lambda=kQ/I$ is biserial (respectively a string algebra, or gentle algebra), then for every $x \in Q_0$ and each $\gamma \in Q_1$ the quotient algebras $\Lambda/\langle e_x \rangle$ and $\Lambda/\langle \gamma \rangle$ are again biserial (respectively string, or gentle). 
Moreover, an arbitrary quotient of a (special) biserial algebra is again (special) biserial.
For $\Lambda = kQ/I$, a \emph{string} in $\Lambda$ is a word $w =\gamma_k^{\epsilon_k}\cdots\gamma_1^{\epsilon_1}$ with letters in $Q_1$ and $\epsilon_i \in \{\pm 1\}$, for all $1 \leq i \leq k$, such that

\begin{enumerate}
\item[(S1)] $s(\gamma_{i+1}^{\epsilon_{i+1}})=e(\gamma_i^{\epsilon_i})$ and $ \gamma_{i+1}^{\epsilon_{i+1}} \neq \gamma_i^{-\epsilon_i}$, for all $1 \leq i \leq k-1$;
\item[(S2)] Neither $w$, nor $w^{-1} := \gamma_1^{-\epsilon_1}\cdots\gamma_k^{-\epsilon_k}$, contain a subpath in $I$.
\end{enumerate}

A string $v$ in $(Q,I)$ is \emph{serial} if either $v$ or $v^{-1}$ is a direct path in $Q$. Namely, $v=\gamma_k \cdots \gamma_2 \gamma_1$ or $v^{-1}=\gamma_k \cdots \gamma_2 \gamma_1$, for some arrows $\gamma_i$ in $Q_1$.
For a string $w=\gamma_k^{\epsilon_k}\cdots\gamma_1^{\epsilon_1}$, we say it starts at $s(w)=s(\gamma_1^{\epsilon_1})$, ends at $e(w)=e(\gamma_k^{\epsilon_k})$, and is of \emph{length} $l(w):=k$. Moreover, a zero-length string, denoted by $e_x$, is associated to every $x \in Q_0$. 
Suppose $\Str(\Lambda)$ is the set of all equivalence classes of strings in $\Lambda$, where for each string $w$ in $\Lambda$ the equivalence class consists of $w$ and $w^{-1}$ (i.e. set $w \sim w^{-1}$).
A string $w$ is called a \emph{band} if $l(w)>0$ and $w^m$ is a string for each $m \in \mathbb{Z}_{\geq 1}$, but $w$ itself is not a power of a string of strictly smaller length, where each band is considered up to all cyclic permutations of it. 
For a vertex $x$ of $\Lambda$, we say $w$ in $\Str(\Lambda)$ \emph{visits $x$} if it is supported by $x$. Moreover, $w$ \emph{passes through $x$} provided that there exists a non-trivial factorization of $w$ at $x$. That is, there exist $w_1, w_2 \in \Str(\Lambda)$ with $s(w_2)=x=e(w_1)$, such that $l(w_1), l(w_2) >0$ and $w=w_2w_1$. For $\alpha \in Q_1$, we say $w$ is \emph{supported} by $\alpha$ if the string $w$ contains $\alpha$ or $\alpha^{-1}$ as a letter.

Let $G_Q$ denote the underlying graph of $Q$. Then, every string $w = \gamma_d^{\epsilon_d}\cdots \gamma_1^{\epsilon_1}$ induces a walk \xymatrix{x_{d+1} \ar@{-}^{\gamma_d}[r] & x_d \ar@{-}^{\gamma_{d-1}}[r] & \cdots & x_1 \ar@{-}_{\gamma_1}[l]} in $G_Q$, where a vertex or an edge may occur multiple times.
The representation $M(w) := ((V_x)_{x \in Q_0}, (\varphi_\alpha)_{\alpha\in Q_1})$ of $(Q,I)$ associated to $w$ has an explicit construction as follows: Put a copy of $k$ at each vertex $x_i$ of the walk induced by $w$. This step gives the vector spaces $\{V_x\}_{x \in Q_0}$, where $V_{x} \simeq k^{n_x}$ and $n_x$ is the number of times $w$ visits $x$.
To specify the linear maps of the representation $M(w)$ between the two copies of $k$ associated to $s(\gamma_i^{\epsilon_i})$ and $e(\gamma_i^{\epsilon_i})$, put the identity in the direction of $\gamma_i$. Namely, this identity map is from the basis vector of $s(\gamma_i)$ to that of $e(\gamma_i)$ if $\epsilon_i=1$, and it goes from the basis vector of $e(\gamma_i)$ to that of $s(\gamma_i)$, if $\epsilon_i=-1$.
If $\Lambda=kQ/I$, the \emph{string module} associated to $w$ is an indecomposable $\Lambda$-module given in terms of the representation $M(w)$. Note that for every string $w$, there is an isomorphism of modules (representations) $M(w) \simeq M(w^{-1})$.

To every band $v \in \Str(\Lambda)$, in addition to the string module $M(v)$, there exists another type of indecomposable $\Lambda$-modules associated to $v$ which are called \emph{band modules}. For the description of band modules, as well as the morphisms between the string and band modules over string algebras, we refer to \cite{BR}, \cite{Kr} and \cite{WW}.

\subsection{$\tau$-tilting (in)finiteness}

Let $M$ be a finitely generated $\Lambda$-module. We say that $M$ is \emph{basic} if no indecomposable module appears more than once in the Krull-Schmidt decomposition of $M$. We denote by $|M|$ the number of non-isomorphic indecomposable direct summands of $M$. Also, $M$ is \emph{rigid} if $\operatorname{Ext}^1_{\Lambda}(M,M)=0$. Let $\rigid (\Lambda)$ denote the set of isomorphism classes of all basic rigid modules in $\modu \Lambda$. 
Moreover, by $\irigid (\Lambda)$ we denote the set of indecomposable modules in $\rigid (\Lambda)$.
Similarly, $M$ is said to be \emph{$\tau$-rigid} if $\Hom_{\Lambda}(M,\tau_{\Lambda} M)=0$, where $\tau_{\Lambda}$ denotes the Auslander-Reiten translation in $\modu \Lambda$. Provided there is no confusion, we simply use $\tau$ to denote the Auslander-Reiten translation. By $\tau \trig(\Lambda)$ and $\mathtt{i}\tau \trig(\Lambda)$ we respectively denote the set of isomorphism classes of basic $\tau$-rigid modules and the indecomposable $\tau$-rigid modules.
A rigid module $X$ is called \emph{tilting} if $\pd_{\Lambda}(X)\leq 1$ and $|X|=|\Lambda|$, where $\pd_{\Lambda}(X)$ denotes the projective dimension of $X$.
Analogously, a $\tau$-rigid module $M$ is \emph{$\tau$-tilting} if $|M|=|\Lambda|$. More generally, $M$ is called \emph{support $\tau$-tilting} if $M$ is $\tau$-tilting over $A/\langle e \rangle$, where $e$ is an idempotent in $A$. 
By $\tilt(\Lambda)$ and $\tau \ttilt(\Lambda)$ we respectively denote the set of all isomorphism classes of basic tilting and $\tau$-tilting modules in $\modu \Lambda$.
Moreover, $s\tau \ttilt(\Lambda)$ denotes the set of isomorphism classes of all basic support $\tau$-tilting modules in $\modu \Lambda$.

$\tau$-tilting theory, introduced by Adachi, Iyama and Reiten \cite{AIR}, has been a modern setup in representation theory of associative algebras where many rich ideas from cluster algebras and classical tilting theory meet. Through this new setting, the authors address the deficiency of classical tilting theory with respect to the mutation of tilting modules. In \cite{AIR}, the notion of mutation of clusters is conceptualized in terms of mutation of (support) $\tau$-tilting modules.

Given an algebra $\Lambda$, it is \emph{a priori} a hard problem to decide whether or not the set of (support) $\tau$-tilting modules is finite. Since these modules form the main ingredient of $\tau$-tilting theory, finding explicit necessary and sufficient conditions such that an algebra has $|\tau \ttilt(\Lambda)| < \infty$ is monumental. This has spurred a lot of research in this direction, among which the elegant ``brick-$\tau$-rigid correspondence"  appearing in \cite{DIJ} has proved to be very useful. 
Some important characterizations of $\tau$-tilting finite algebras are recalled in the rest of this subsection.

Recall that a $\Lambda$-module $Y$ is called a \emph{brick} if $\End_{\Lambda}(Y)$ is a division algebra. That is, any non-zero endomorphism of $Y$ is invertible. As in Section \ref{Section:Introduction}, by $\Brick(\Lambda)$ and $\brick(\Lambda)$ we respectively denote the set of isomorphism classes of bricks in $\Modu \Lambda$ and $\modu \Lambda$. If the field $k$ is algebraically closed, then $Y$ belongs to $\brick (\Lambda)$ if and only if $\End_{\Lambda}(Y)\simeq k$.
Such modules are sometimes called \emph{Schur representations}, particularly when they are studied from the algebro-geometric viewpoint, such as in \cite{CKW}.
An algebra $\Lambda$ is called \emph{brick-finite} provided $|\Brick(\Lambda)| < \infty$. Meanwhile, we warn the reader that those algebras called ``Schur representation-finite" in \cite{CKW} are not known to be necessarily brick-finite (for further details on this difference, see Subsection \ref{Subsection:Schemes and varieties of representations}, as well as \cite[Subsection 1.3]{Mo2}).

One of our main goals in this paper is to establish a relationship between certain modules in $\Brick(\Lambda)\setminus \brick(\Lambda)$ and those in $\brick(\Lambda)$.
In this regard, the following result of Sentieri \cite{Se} is of interest.

\begin{theorem}[\cite{Se}]\label{Thm:Sentieri}
An algebra $\Lambda$ is brick-finite if and only if every brick in $\Modu \Lambda$ is finite dimensional. 
\end{theorem}

We now list some of the fundamental results on $\tau$-tilting finiteness of algebras. Recall that a subcategory $\mathcal{T}$ of $\modu \Lambda$ is a \emph{torsion class} if it is closed under quotients and extensions. Let $\tors(\Lambda)$ denote the set of all torsion classes in $\modu \Lambda$.
For $M$ in $\modu \Lambda$, let $\Fac(M)$ denote the subcategory of $\modu \Lambda$ consisting of all those modules that are quotients of some finite direct sum of copies of $M$.
It is known that $\mathcal{T}$ in $\tors(\Lambda)$ is \emph{functorially finite} provided $\mathcal{T}=\Fac(M)$, for some $M$ in $\modu \Lambda$. By $\ftors(\Lambda)$ we denote the subset of $\tors(\Lambda)$ consisting of functorially finite torsion classes. The following important result relates the finiteness of the notions introduced so far. In particular, it states that an algebra is brick-finite if and only if it is $\tau$-tilting finite.

\begin{theorem}(\cite{AIR, DIJ})\label{Thm:tau-tilting finiteness equivalences}
For an algebra $\Lambda$, the following are equivalent:
\begin{enumerate}
\item  $\Lambda$ is $\tau$-tilting finite;
\item  $s\tau \ttilt(\Lambda)$ is finite;
\item $\tau\trig(\Lambda)$ is finite;
\item $\brick(\Lambda)$ is finite;
\item  $\tors(\Lambda)=\ftors(\Lambda)$.
\end{enumerate}
\end{theorem}

\subsection{Minimal brick-infinite algebras}\label{Subsection:Minimal brick-infinite algebras}
Here we collect some of our main results from \cite{Mo1, Mo2}, as well as \cite{MP}, which are used in this paper. We begin with a useful observation that is freely used in our reductive arguments. In particular, we recall that each epimorphism of algebras $\psi: \Lambda_1 \rightarrow \Lambda_2$ induces an exact functorial full embedding $\widetilde{\psi}: \modu\Lambda_2 \rightarrow \modu\Lambda_1$. Particularly, we get $\ind(\Lambda_2) \subseteq \ind(\Lambda_1)$ and also $\brick(\Lambda_2) \subseteq \brick(\Lambda_1)$.
This implies that if $\Lambda_2$ is rep-infinite (respectively, brick-infinite) then so is $\Lambda_1$.
Thus, by Theorem \ref{Thm:tau-tilting finiteness equivalences}, $\tau$-tilting finiteness is preserved under taking quotients.

Recall that an algebra $\Lambda$ is \emph{minimal representation-infinite} (or min-rep-infinite, for short) if $\Lambda$ is rep-infinite and any proper quotient of $\Lambda$ is representation-finite. Following our notations in \cite{Mo2}, by 
$\Mri({\mathfrak{F}_{\sB}})$ we respectively denote the family of min-rep-infinite special biserial algebras and $\Mri({\mathfrak{F}_{\nD}})$ denotes the family of non-distributive min-rep-infinite algebras. 
Before we summarize the relevant results on the brick-(in)finiteness of these algebras, let us recall that in \cite{Mo1}, the following bound quivers are called \emph{generalized barbell}: 

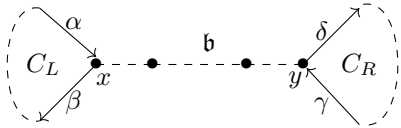
\begin{center}
\begin{tikzpicture}

 \draw [->] (1.25,0.75) --(2,0.1);
    \node at (1.7,0.55) {$\alpha$};
 \draw [<-] (1.25,-0.75) --(2,0);
    \node at (1.7,-0.5) {$\beta$};
  \draw [dashed] (1.25,0.75) to [bend right=100] (1.25,-0.75);
   \node at (1.3,0) {$C_L$};
    \node at (2,0) {$\bullet $};
    \node at (2.1,-0.2) {$x$};
 \draw [dashed] (2,0) --(2.75,0);
 \node at (2.75,0) {$\bullet$};
 \draw [dashed] (2.75,0) --(4,0);
 \node at (4,0) {$\bullet$};
 \draw [dashed] (4,0) --(4.75,0);
 \node at (3.5,0.3) {$\mathfrak{b}$};
 
 \node at (4.75,0) {$\bullet$};
 \draw [<-] (5.5,0.75) --(4.75,0);
 \node at (5,0.45) {$\delta$};
 \draw [->] (5.50,-0.8) --(4.8,-0.05);
    \node at (5,-0.55) {$\gamma$};
  \draw [dashed] (5.55,0.8) to [bend left=100] (5.55,-0.8);
   \node at (5.5,0) {$C_R$};
   \node at (4.65,-0.2) {$y$};
\end{tikzpicture}
\end{center}
where $ I= \langle \beta \alpha , \delta \gamma \rangle$, $C_L= \alpha \cdots \beta$ and $C_R=\gamma \cdots \delta$ are cyclic strings each of which having no repeated vertex (except for the endpoints) and where $C_L, C_R$ have no common vertex, except for possibly the case where $\mathfrak{b}$ is of zero length (which implies $x=y$). Moreover, $\mathfrak{b}$ (respectively $C_L$ and $C_R$) can have any length (respectively any positive length) and arbitrary orientation of their arrows, provided $C_RC_L$ is not a serial string in $(Q,I)$. The latter case occurs exactly when $\mathfrak{b}$ is of length zero and both $C_L$ and $C_R$ are serial strings. Observe that the assumption that $C_RC_L$ is not a serial string guarantees that $kQ/I$ is a finite dimensional algebra.
We note that generalized barbell quivers are a slight generalization of ``barbell" quivers introduced by Ringel \cite{Ri1}, where he always assume the bar $\mathfrak{b}$ is of positive length.

The following theorem summarizes some of our earlier results on the study of $\tau$-tilting finiteness. To make them more congruent with the scope of this paper, below we state them in terms of bricks.

\begin{theorem}[\cite{Mo2}]\label{tau-inf min-rep-infinite classification}
With the same notations as above, the following hold:

\begin{enumerate}
    \item\label{tau-inf min-rep-infinite sp.bis} If $\Lambda$ belongs to $\Mri({\mathfrak{F}_{\sB}})$, then $\Lambda$ is brick-infinite if and only if $(Q,I)$ is hereditary of type $\widetilde{A}_n$ or $(Q,I)$ is generalized barbell.
    
    \item If $\Lambda$ belongs to $\Mri({\mathfrak{F}_{\nD}})$, then $\Lambda$ is brick-infinite if and only if $Q$ has a sink.
    
\end{enumerate}
\end{theorem}
We remark that $\Mri({\mathfrak{F}_{\sB}})$ and $\Mri({\mathfrak{F}_{\nD}})$ consist of only tame algebras and either of these two families contains both brick-finite and brick-infinite algebras (see \cite{Mo1, Mo2} for full classifications).

We also recall that an algebra $\Lambda$ is said to be \emph{minimal $\tau$-tilting infinite} if $\Lambda$ is $\tau$-tilting infinite but every proper quotient of $\Lambda$ is $\tau$-tilting finite. From Theorem \ref{Thm:tau-tilting finiteness equivalences}, it is immediate that minimal $\tau$-tilting infinite algebras are the same as minimal brick-infinite algebras. Here we only list some of the main properties of these algebras and for more details we refer to \cite{MP}.
Recall that $\Lambda$ is called \emph{central} provided its center is the ground field $k$.

\begin{theorem}\label{min-tau-inf properties}
Let $\Lambda=kQ/I$ be a minimal brick-infinite algebra. Then, 
\begin{enumerate}
    \item $\Lambda$ is central and admits no projective-injective module. Moreover, $(Q,I)$ has no node.
    \item Almost every $\tau$-rigid $\Lambda$-module is faithful, and therefore is partial tilting.
    \item $\Lambda$ is minimal tilting infinite (i.e. $\tilt(\Lambda)$ is an infinite set but $\tilt(A/J)$ is finite, for each non-zero ideal $J$ in $\Lambda$).
\end{enumerate}
\end{theorem}

To highlight some fundamental differences between these modern and classical notions of minimality, we remark that min-rep-infinite algebras are not necessarily central and their bound quivers can have several nodes. Furthermore, note that although $\tau$-tilting finiteness is preserved under algebraic quotients, there exists tilting-finite algebra $\Lambda$ such that $\Lambda/J$ is tilting-infinite, for an ideal $J$ in $\Lambda$.

\medskip

\subsection{Schemes and varieties of representations}\label{Subsection:Schemes and varieties of representations}
In this subsection we collect some basic tools used in this paper which allow us to move between the algebraic and geometric sides of our problem.
In particular, for algebra $\Lambda$ and a dimension vector $\underline{d}$ in $\mathbb{Z}^{|\Lambda|}_{\geq 0}$, let $\rep(\Lambda, \underline{d})$ denote the affine (not necessarily irreducible) variety parametrizing the modules in $\modu (\Lambda, \underline{d})$. Here, $\modu (\Lambda, \underline{d})$ denotes the subcategory of $\modu \Lambda$ consisting of all modules of dimension vector $\underline{d}$.

Under the action of $\GL(\underline{d})$ via conjugation, $\rep (\Lambda, \underline{d})$ can be viewed as a scheme, as well as an affine variety, where the orbits of this action are in bijection with the isomorphism classes of modules in $\modu(\Lambda, \underline{d})$. Through this conceptual dictionary, we study some geometric properties of representations of the bound quiver $(Q,I)$, where $\Lambda=kQ/I$ is an admissible presentation of $\Lambda$.
For $M$ in $\modu(\Lambda, \underline{d})$, by $\mathcal{O}_M$ we denote the $\GL(\underline{d})$-orbit of $M$, when it is viewed as a point in $\rep (\Lambda, \underline{d})$. 
If $\rep (\Lambda, \underline{d})$ is viewed as the $k$-points of a corresponding scheme, it is known that $\mathcal{O}_M$ is open in this scheme if and only if $M$ is rigid. However, if we consider $\rep (\Lambda, \underline{d})$ as a variety, there could be non-rigid modules $N$ such that $\mathcal{O}_N$ is open. 
Although both of these geometric structures are rich and come with powerful tools, we mostly treat $\rep (\Lambda, \underline{d})$ as an affine variety. 
When there is no risk of confusion, $\modu(\Lambda, \underline{d})$ is referred to as a variety to reflect the geometric structure that comes from $\rep (\Lambda, \underline{d})$.

Let $\ind(\Lambda, \underline{d})$ and $\brick(\Lambda, \underline{d})$ respectively denote the set of all indecomposable modules and bricks in $\modu (\Lambda,\underline{d})$. It is known that $\brick(\Lambda, \underline{d})$ is an open subset of $\modu (\Lambda,\underline{d})$. Let $\Irr(\Lambda,\underline{d})$ be the set of all irreducible components of $\modu(\Lambda, \underline{d})$, and by $\Irr(\Lambda)$ we denote the union of all $\Irr(\Lambda,\underline{d})$, where $\underline{d}$ is an arbitrary dimension vector. A component $\mathcal{Z} \in \Irr(\Lambda)$ is called \emph{indecomposable} provided it contains a non-empty open subset $U$ which consists of indecomposable representations. 
In \cite{CBS}, the authors prove a geometric analogue of the Krull-Schmidt decomposition for irreducible components, which highlights the role of indecomposable components among all irreducible ones.

For each $\mathcal{Z}$ in $\Irr(\Lambda)$, the algebraic properties of the modules in $\mathcal{Z}$ capture important information on the geometry of $\mathcal{Z}$, and vice versa.
Motivated by this interaction, Chindris, Kinser and Weyman \cite{CKW} have recently adopted a geometric approach to generalize the notion of representation-finiteness, primarily based on the properties of irreducible components.
In particular, $\Lambda$ is said to have \emph{dense orbit property} provided every $\mathcal{Z}$ in $\Irr(\Lambda)$ contains a dense orbit. By some simple geometric considerations, one can show that every rep-finite algebra has the dense orbit property. In \cite{CKW}, the authors show that the new notion is novel and construct explicit rep-infinite algebras which have the dense orbit property. Furthermore, they prove that a string algebra (and more generally, each special biserial algebra) is rep-finite if and only if it has the dense orbit property.

Adopting this algebro-geometric approach, we say that $\Lambda$ is \emph{brick-discrete} if for each $d \in \mathbb{Z}_{\geq 0}$, there are only finitely many (isomorphism classes of) bricks of dimension $d$. This is equivalent to the fact that for each $\mathcal{Z}$ in $\Irr(\Lambda)$, if $M$ belongs to $\brick(\mathcal{Z})$, then $\mathcal{Z}=\overline{\mathcal{O}}_M$. Here, $\overline{\mathcal{O}}_M$ denotes the orbit closure of $\mathcal{O}_M$.
We remark that brick-discrete algebras have been treated in \cite{CKW}, where the authors introduced them under the name ``Schur representation-finite" algebras. To avoid confusion between brick-finite and Schur representation-finite algebras, we use our new terminology and call the latter type brick-discrete.

In \cite{CKW}, it is shown that if $\Lambda$ has the dense orbit property, then it is brick-discrete, but the converse does not hold in general. So, brick-discreteness was considered as a generalization of the dense orbit property, and hence a generalization of rep-finiteness. 
We observe that every brick-finite algebra is brick-discrete.
In contrast, in general it is not known whether brick-discrete algebras are necessarily brick-finite.
In fact, this is the content of the following conjecture, which is a precursor of Conjecture \ref{Introduction: Conjecture}.
\begin{conjecture}[\cite{Mo2}]\label{My Old Conjecture}
Let $\Lambda$ be an algebra over an algebraically closed field $k$. The following are equivalent:
\begin{enumerate}
    \item $\Lambda$ is $\tau$-tilting finite;
    \item $\Lambda$ is brick-discrete;
\end{enumerate}
Equivalently, $\Lambda$ is brick-infinite if and only if there exists a family $\{M_{\lambda}\}_{\lambda \in k^*}$ of bricks of the same length.
\end{conjecture}

The above conjecture first appeared in the arXiv version of \cite{Mo2}, where the first-named author proposed an algebro-geoemtric realization of $\tau$-tilting finiteness. Moreover, it is verified for all algebras treated in that paper. We also remark that the numerical implication of the above conjecture was later stated in \cite{STV}.

\medskip

For a special biserial algebra $\Lambda$, an irreducible component $\mathcal{Z}$ in $\Irr(\Lambda)$ is called a \emph{string component} if it contains a string module $M$ such that $\mathcal{O}_M$ is dense in $\mathcal{Z}$. That being the case, we get $\mathcal{Z}=\overline{\mathcal{O}}_M$, which implies that $\mathcal{Z}$ can be specified by the isomorphism class of the string module $M$. 
In contrast, $\mathcal{Z}$ is a \emph{band component} provided $\mathcal{Z}$ contains a family $\{M_{\lambda}\}_{\lambda\in k^*}$ of band modules such that $\bigcup_{\lambda\in \mathcal{C}}\mathcal{O}_{M_{\lambda}}$ is dense in $\mathcal{Z}$. In this case, $\mathcal{Z}=\overline{\bigcup_{\lambda\in \mathcal{C}}\mathcal{O}_{M_{\lambda}}}$. Hence, a band component is determined by the band that gives rise to the one-parameter family $\{M_{\lambda}\}_{\lambda\in k^*}$.
Provided $\Lambda$ is a string algebra, $\Irr(\Lambda)$ consists only of string and band components.

\section{Minimal brick-infinite biserial algebras}\label{Section:Minimal brick-infinite biserial algebras}

It is well-known that a minimal representation-infinite algebra which is special biserial must be a string algebra. Recently, Ringel \cite{Ri1} gave a full classification of these algebras and one can observe that in fact every biserial min-rep-infinite algebra falls into Ringel's classification (for further details, see \cite{Mo2}).
In this section, we give an analogous classification result and fully describe the bound quivers of those biserial algebras which are minimal brick-infinite.
In particular, we show that any minimal brick-infinite biserial algebra is gentle and falls into exactly one of the two types described in Theorem \ref{Thm: Classification of min-brick-infinite biserials}.

We first recall some notations and results from \cite{CB+}. A quiver $Q$ is said to be \emph{biserial} if for any vertex $x$ in $Q$, there are at most two arrows starting at $x$, and at most $2$ arrows ending at $x$. It is clear that the quiver of any biserial algebra has to be biserial. A \emph{bisection} $(\sigma, \tau)$ of a biserial quiver $Q$ is the data of two functions $\sigma, \tau: Q_1 \to \{\pm 1\}$ such that if $\alpha, \beta$ are two distinct arrows starting (resp. ending) at $x$, then $\sigma(\alpha) \ne \sigma(\beta)$ (resp. $\tau(\alpha) \ne \tau(\beta)$). Given a biserial quiver and a bisection of it, a \emph{good path} is any path $\alpha_r \cdots \alpha_1$ such that for $1 \le i \le r-1$, we have that $\tau(\alpha_i) = \sigma(\alpha_{i+1})$. Trivial paths are declared to be good. A path that is not good is said to be \emph{bad}. Bad paths of length two will play an important role due to the following result.

Observe that if $\Lambda=kQ/I$ is such that there exist multiple arrows between two fixed vertices of $Q$, then $\Lambda$ is minimal brick-infinite if and only if $Q$ is the Kronecker quiver and $I=0$.
Hence, for simplicity of the assertions, in the rest of this section we exclude the situation where the quiver of algebra has multiple arrows.

\begin{theorem}[\cite{CB+}] \label{DescriptionBiserial}
Let $A$ be a biserial algebra with quiver $Q$ having no multiple arrows. There exists a bisection $(\sigma, \tau)$ of $Q$ such that $A \cong kQ/I$, and for each bad path $\beta\alpha$ of length two $I$ contains an element $R_{\beta\alpha}$ of one of the following types: 
\begin{enumerate}[$(1)$]
    \item $R_{\beta\alpha} = \beta\alpha$ or
    \item there is a path $p$ parallel to $\beta$, which neither starts nor ends with $\beta$ such that $p\alpha$ is good and such that $R_{\beta\alpha} = \beta\alpha - \lambda_{\beta\alpha}p\alpha$ for some non-zero scalar $\lambda_{\beta\alpha}$.
\end{enumerate} 
Conversely, if $Q$ is a biserial quiver with no multiple arrows with a bisection $(\sigma, \tau)$, and $I$ is an admissible ideal of $kQ$ that contains all of the above elements $R_{\beta\alpha}$, then $kQ/I$ is biserial.
\end{theorem}

Since the opposite algebra to a biserial algebra is also biserial, the dual of the above theorem is also valid as follows.

\begin{theorem} \label{DescriptionBiserialDual}
Let $A$ be a biserial algebra with quiver $Q$ having no multiple arrows. There exists a bisection $(\sigma, \tau)$ of $Q$ such that $A \cong kQ/I$, and for each bad path $\beta\alpha$ of length two $I$ contains an element $R'_{\beta\alpha}$ of one of the following types:
\begin{enumerate}[$(1)$]
    \item $R'_{\beta\alpha} = \beta\alpha$ or
    \item there is a path $p'$ parallel to $\alpha$, which neither starts nor ends with $\alpha$ such that $\beta p'$ is good and such that $R'_{\beta\alpha} = \beta\alpha - \lambda'_{\beta\alpha}\beta p'$ for some non-zero scalar $\lambda'_{\beta\alpha}$.
\end{enumerate} 
Conversely, if $Q$ is a biserial quiver with no multiple arrows with a bisection $(\sigma, \tau)$, and $I$ is an admissible ideal of $kQ$ that contains all of the above elements $R'_{\beta\alpha}$, then $kQ/I$ is biserial.
\end{theorem}

\medskip

The next result plays an important role in our arguments below. This proposition appears as Lemma 2.3 in \cite{Ku}, but it is originally due to Vila-Freyer.

\begin{proposition}[Vila-Freyer] \label{BadPathSocle}
Consider the setting in the above theorems. Then, any bad path $\beta\alpha$ is such that for any arrow $\gamma$, we have 
\begin{enumerate}[$(1)$]
    \item $\gamma\beta\alpha \in I$ in the setting of Theorem \ref{DescriptionBiserial}.
    \item $\beta\alpha\gamma \in I$ in the setting of Theorem \ref{DescriptionBiserialDual}.
\end{enumerate}
\end{proposition}

It is important to note that in the setting of Theorem \ref{DescriptionBiserial}, although we always have $\gamma \beta\alpha \in I$, there may exist an arrow $\delta$ such that $\beta \alpha \delta \notin I$. A similar observation holds for Theorem \ref{DescriptionBiserialDual}. The following example further explains this phenomenon.

\begin{example}
Let $Q$ be the following biserial quiver
$$
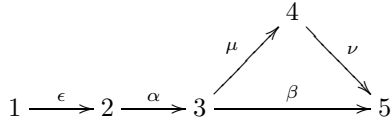
\xymatrix{&& & 4 \ar[dr]^\nu & \\ 1 \ar[r]^\epsilon & 2 \ar[r]^\alpha & 3 \ar[rr]^\beta \ar[ur]^\mu && 5 }$$
with bisection such that the values of $\sigma, \tau$ are always positive, except that $\sigma(\beta) = \tau(\beta) = -1$. There is a unique bad path of length $2$, namely $\beta\alpha$. The algebra $kQ/\langle \beta\alpha - \nu\mu\alpha\rangle$ is biserial by Theorem \ref{DescriptionBiserial}. However, the bad path $\beta\alpha\epsilon$ is not in $I$. On the other hand, one can change the generator $\beta$ by $\beta' = \beta - \nu\mu$. Using the same bisection, now we get that the bad path $\beta'\alpha$  belongs to $I$. Hence, $A = kQ/I$ is special biserial. On the other hand, let $Q'$ be the quiver obtained from $Q$ by adding a vertex $6$ and an arrow $\gamma: 5 \to 6$. Consider the bisection of $Q'$ obtained by extending the bisection of $Q$ with $\sigma(\gamma) = \tau(\gamma)=1$. Consider the algebra $A' = kQ' / \langle \beta\alpha - \nu\mu\alpha, \gamma\beta\rangle$. Then the algebra $A'$ is biserial but no longer special biserial. With the change of generator as above, the relations $\beta\alpha - \nu\mu\alpha$ and $\gamma\beta$ respectively become $\beta'\alpha$ and  $\gamma\beta' + \gamma\nu\mu$. The latter is a relation as in Theorem \ref{DescriptionBiserialDual}.
\end{example}

The rest of this section is dedicated to proving the following theorem.
\begin{theorem}\label{Thm: Classification of min-brick-infinite biserials}
Let $A$ be a minimal brick-infinite biserial algebra. Then, $A$ is a generalized barbell algebra, or $A=k\widetilde{A}_m$, for some $m\in \mathbb{Z}$.
\end{theorem}

To prove the above theorem, in the remainder of the section we assume $A$ is minimal brick-infinite and is not hereditary and show that $(Q,I)$ must be a generalized barbell quiver. Hence, below $A=kQ/I$ always denotes a minimal brick-infinite biserial algebra which is not hereditary. In particular, the quiver $Q$ does not have multiple arrows (between any pair of vertices).
We fix a bisection $(\sigma, \tau)$ of $Q$ and assume that $A \cong kQ/I$ where $I$ is as Theorem \ref{DescriptionBiserial}. We may further assume that the bisection and ideal as in Theorem \ref{DescriptionBiserial} (and similarly, in Theorem \ref{DescriptionBiserialDual}) are chosen in such a way that, up to isomorphism, the number of relations of type $(2)$ in that theorem is minimal. 
We also note that for every minimal brick-infinite algebra $A$, and for each non-zero element $r$ in $A$, only finitely many non-isomorphic bricks in $\modu A$ are anniliated by $r$. This is because if we assume otherwise, the quotient algebra $A/\langle r \rangle$ becomes brick-infinite and this contradicts the minimality assumption on $A$. We freely use this property of minimal brick-infinite algebras, including in the proof of the following lemma.

\begin{lemma}\label{Lemma: particular relations}
In the above setting, if $A$ is not special biserial, there has to be one relation of type $(2)$ from Theorem \ref{DescriptionBiserial} such that $\beta$ is not a loop.
\end{lemma}

\begin{proof}

For the sake of contradiction, assume otherwise. That is, for any relation $R_{\beta\alpha} = \beta\alpha - \lambda p \alpha$ with non-zero $\lambda$, where $p\alpha$ is a good path and $\beta\alpha$ is a bad path of length two passing through $x$, we must have $\beta$ is a loop.
Since $p$ neither starts nor ends with $\beta$, there are arrows $u, v$, different from $\beta$, such that $p = vp'u$ for some good path $p'$. If $u = v$ is a loop, then $Q$ has a single vertex with two loops, $A$ is local and hence brick-finite, a contradiction. Hence, $u \ne v$ are not loops.

If $\alpha = \beta$, we get $\beta^2 = \lambda p \beta$. In this case, it is clear that $uv$ and $\beta^2$ are bad paths and $\beta$ does not appear in any other bad path of length two. By our assumption, we have $uv \in I$. We set $\beta' = \beta - \lambda p$. Assume first that $(\beta')^2 \in I$. Consider the change of generators which changes only $\beta + I$ to $\beta' + I$, so that the arrow originally representing $\beta + I$ now represents $\beta' + I$. Observe that after this change of generators, the bisection is preserved. Since both bad paths at $x$ are now represented by elements in $I$, up to isomorphism, we have that $A \cong kQ/J$ where the number of relations of type $(2)$ in $J$ from Theorem \ref{DescriptionBiserial} has decreased by one, contradicting our running minimality assumption on the number of such relations. Hence, we may assume $(\beta')^2 \not \in I$. We see that $\gamma(\beta')^2 \in I$ for $\gamma = u$ and $\gamma = \beta'$ by Proposition \ref{BadPathSocle}. Therefore, any module not annihilated by $(\beta')^2$ has a submodule isomorphic to the simple at $x$. On the other hand, it can be shown by direct computations that $(\beta')^2\gamma \in I$, for any arrow $\gamma$. This yields that any module not annihilated by $(\beta')^2$ has a quotient isomorphic to the simple at $x$. Therefore, for any brick $M$, we have that $M$ is annihilated by $(\beta')^2$. This gives the desired contradiction, because over the minimal brick-infinite algebra $A$, at most finitely many bricks are annihilated by a fixed non-zero element.

Assume now that $\alpha$ and $\beta$ are distinct. Therefore, $\alpha$ is not a loop and $\alpha = v$. The bad paths at $x$ are $\beta\alpha$ and $u\beta$. Because $u$ is not a loop and $u\beta$ is bad, we get $u\beta \in I$. As in the previous case, we set $\beta' = \beta - \lambda p$ and make a similar argument. If $u\beta' \in I$, then we proceed with a change of generator for the loop $\beta$ (so $\beta + I$ is replaced by $\beta' + I$) and we get another presentation of $A$ with a smaller number of relations of type $(2)$, leading to a contradiction. If $u\beta'$ is not in $I$, then $u \beta'\gamma \in I$ for $\gamma = \beta'$ and $\gamma = \alpha$ by direct computations. Also, it follows from Proposition \ref{BadPathSocle} that $\gamma\beta\alpha \in I$, for any arrow $\gamma$. Now, any module not annihilated by $\beta\alpha$ has a submodule isomorphic to the simple at $x$ and any module not annihilated by $u\beta'$ has a quotient isomorphic to the simple at $x$. We get a contradiction as in the previous case.
\end{proof}

From the preceding lemma and our assumptions, it follows that if $A$ is minimal brick-infinite and not special biserial, then $A=kQ/I$ with a bisection of $Q$ such that there exists a bad path $\beta\alpha$, where $\beta$ is not a loop and where $R_{\beta\alpha} = \beta\alpha - \lambda p \alpha \in I$ for some non-zero scalar $\lambda$, and a good path $p\alpha$. We need the following lemmas.

\begin{lemma} \label{MonomialRelationsBiserial}
Assume that $A = kQ/I$ where $I$ contains the generators from Theorem \ref{DescriptionBiserial}. Then, $I$ is generated by these relations, plus possibly some other monomial relations that are all good paths.
\end{lemma}

\begin{proof}
 Let $I'$ be the ideal containing the generators from Theorem \ref{DescriptionBiserial}, and assume that $I' \subset I$ is a proper inclusion. Take $R$ in $I \setminus I'$  as a relation that starts at vertex $c$ and ends at vertex $d$ (where we may have $c=d$).  
Observe that modulo $I'$, every bad path is a (possibly zero) scalar multiple of a good path.
We may assume from the beginning that $R$ is a linear combination of good paths with no summand that appears in $I'$. Note that $R$ can have more than $2$ terms.

We first treat the case where all of the paths in the summands of $R$ start with the same arrow. In this case, we can find a good path $p$ and write the relation as $R = \lambda p + xp$, for some scalar $\lambda \ne 0$ and $x$ a linear combination of paths of length at least one. This implies that $p+I$ lies in arbitrary large powers of the radical of $A$. Since the radical is nilpotent, we conclude that $p \in I$, and our assumption on the summands of $R$ implies that $R$ is monomial. 

Now, let us assume that the paths occurring in $R$ start with two different arrows. Consequently, there are paths $p$ and $q$ starting at $c$ with two different arrows such that $R = vq + up$ where $vq \notin I$ and $up \notin I$. We remark that $u$ (and similarly $v$) is not necessarily a single path, but it is a linear combination of paths. 
Since $A$ is biserial, the radical of the projective module $Ae_c$ is a sum of two uniserial modules $U_1$ and $U_2$, whose intersection has to be the simple module at vertex $d$.
Hence, $Ae_c$ embeds into the injective module $D(e_dA)$ at $d$. Using dual argument, the right projective module $e_dA$ embeds into the right injective module $D(Ae_c)$. By dimension count, that means $Ae_c$ is projective-injective, which cannot happen for a minimal brick-infinite algebra, by 
Theorem \ref{min-tau-inf properties}(1).
\end{proof}

In the next lemma, we describe the behaviour of certain parallel paths in the bound quiver of algebras under consideration. 
We reiterate that in the following lemma, we still work with the minimal brick-infinite algebra $A=kQ/I$ under the running assumptions that we fixed in the paragraph following Theorem \ref{Thm: Classification of min-brick-infinite biserials}.

\begin{lemma} \label{LemmaParallel}
Let $p_1, p_2$ be two parallel good paths starting at a vertex $c$ and ending at a vertex $d$. Assume further that $p_i \not \in I$ and no arrow on $p_i$ starts with vertex $d$, for $i=1,2$. Then, there is an arrow of $p_1$ or $p_2$ ending at $c$.
\end{lemma}
\begin{proof}
We assume otherwise and claim that the good paths $p_1$ and $p_2$ are linearly independent as elements in the $e_dAe_d-e_cAe_c$-bimodule $e_dAe_c/\rad(e_dAe_c)$. Here, $\rad(e_dAe_c)$ denotes the radical of the bimodule $e_dAe_c$.  To verify the claim, note that if $p_1$ and $p_2$ are linearly dependent, there are elements $u \in e_cAe_c$ and $v \in e_dAe_d$ such that $p_1 = vp_2u$ or $p_2 = vp_1u$. 

Without loss of generality, we assume $p_1 = vp_2u$. By Lemma \ref{MonomialRelationsBiserial}, $I$ is generated by the generators from Theorem \ref{DescriptionBiserial}, plus possibly some other monomial relations. Note that from the proof of Lemma \ref{MonomialRelationsBiserial}, these monomial relations can all be taken to be good paths.
We have that $p_1$ is a good path and $p_1 \notin I$. Therefore, $u \notin I$ and $u$ can be taken to be a linear combination of good paths from $c$ to $c$. The latter is an immediate consequence of the presentation described in Theorem \ref{DescriptionBiserial}.
Now, since none of $p_1, p_2$ revisit vertex $c$, each non-trivial good path starting at $c$ and returning to $c$ either starts with $p_1$ or with $p_2$. Therefore, we get an expression $p_1 = w_1p_1 + w_2p_2$ where $w_1, w_2 \in A$ and $w_1$ is in the radical. This implies that for each positive integer $n$, we have an expression of the form $p_1 = w_1^np_1 + u_np_2$ for some $u_n \in A$. Note that for $n$ large enough, $w_1^n \in I$. 
Therefore, $p_1 = w_3p_2$, for some $w_3 \in A$. 

Now, consider the projective module $Ae_c$. If there is only one arrow starting at $c$, then the longer of $p_1$ or $p_2$ contains the shorter one as a subpath, which contradicts the assumption of the lemma. Since there are two arrows starting at $c$, we have that ${\rm rad}A e_c$ is a sum of two uniserial modules, and $p_1 = w_3p_2$ implies that these two uniserial modules intersect non-trivially at a simple submodule at $d$. Hence, $Ae_c$ embeds into the injective module $D(e_dA)$ at $d$. 

Dually, arguing on $v$, we similarly get that $p_1 = p_2w_3''$ where $w_3'' \in A$. As above, this yields that $e_dA$ embeds into the injective module $D(Ae_c)$. By dimension count, we get $Ae_c$ is projective-invective, which cannot happen for a minimal brick-infinite algebra, by Theorem \ref{min-tau-inf properties}(1).
\end{proof}

In what comes next, we will use the following setting.

\begin{Setting} \label{setting1}
\begin{enumerate}
\item $A$ is a biserial algebra given by a bound quiver $(Q,I)$ and there is a bisection $(\sigma, \tau)$ of $Q$ such that $I$ is as Theorem \ref{DescriptionBiserial}.
\item $A$ is minimal brick-infinite and not hereditary.
\item The bisection $(\sigma, \tau)$ and ideal $I$ are chosen so that the number of relations of type (2) from Theorem \ref{DescriptionBiserial} is minimal.
\item There is a bad path $\beta\alpha$ of length two, where $\beta$ is not a loop, by Lemma \ref{Lemma: particular relations}. Also, we have $R_{\beta\alpha} = \beta\alpha - \lambda p\alpha \in I$, where $\lambda$ is a non-zero scalar and $p\alpha$ is a good path.
\end{enumerate}
\end{Setting}

\medskip 

\begin{lemma} \label{beta_in_p}
    Using Setting \ref{setting1}, $\beta$ must appear in $p$.
\end{lemma}  

\begin{proof}
Assume otherwise, that is, suppose $\beta: a \to b$ does not appear in the good path $p$. Then we take a starting subpath $p'$ of $p$ (thus starting at $a$) which is minimal with the property of ending at $b$. Applying Lemma \ref{LemmaParallel} to $\beta$ and $p'$ yields that $a$ has to be the ending vertex of an arrow in $p'$. Hence, there is a non-trivial good path $q$ from $a$ to $a$ such that $p' = p''q$ where $p''$ is a good path starting at $a$ and ending at $b$, and it does not go through $\beta$. But this means that $p''$ starts with $p'$, which yields an expression $p''= q'p'$, for a good path $q'$, and therefore $p' = q'p'q$. This is absurd and completes the proof. 
\end{proof}

We recall that a generalized barbell algebra whose bar is of positive length is simply called a barbell algebra (see Subsection \ref{Subsection:Minimal brick-infinite algebras}). In the next definition, we specify a particular type of barbell algebras which are important in the proof of our main theorem.

\begin{definition}
Consider the quiver $Q$ with given bisection and ideal $I$ as in Theorem \ref{DescriptionBiserial} so that $A = kQ/I$. We say that a subquiver $Q'$ is a \emph{good barbell subquiver} for $A$ if 
\begin{enumerate}
    \item $Q'$ is the quiver of a barbell algebra such that the two defining cycles are oriented cycles and the bar is linear (that is, the bar is a path);
    \item any path which is not going through the ``zero relations" of the barbell is a non-zero good path.
\end{enumerate} 
\end{definition}

\medskip 

\begin{construction}\label{Remark: Existence of good barbell subquiver}
We assume that Setting \ref{setting1} holds. Our goal is to show that there is a good barbell subquiver of $Q$. According to Lemma \ref{beta_in_p}, the good path $p$ contains $\beta$, and it is neither the first, nor the last arrow of $p$ by Theorem \ref{DescriptionBiserial}. We construct a good barbell subquiver of $Q$ as follows. 
We encourage the reader to consult Figure \ref{fig: good barbell quiver} while going through the construction below and the arguments following that.
\medskip

First, write the good path $p$ as $p = p_2\beta p_1$ where $\beta$ does not appear in $p_1$. We further write $p_1 = p_{12}p_{11}$ such that $p_{12}$ is the ending subpath of $p_1$ and it is minimal with the property of revisiting itself. Thus, we can write $p_{12}$ as $q_2q_1$ where $q_1$ is a simple oriented cycle from a vertex $a'$ to itself, while $q_2$ is a path from $a'$ to $a$ (if $a = a'$, then $q_2$ is trivial). Similarly, we write $p_2 = p_{22}p_{21}$ where $p_{21}$ is the starting subpath of $p_2$ which is minimal with the property of revisiting itself. Therefore, we can write $p_{21}$ as $q_4q_3$ where $q_4$ is a simple oriented cycle from a vertex $b'$ to itself while $q_3$ is a path from $b$ to $b'$ (if $b = b'$, then $q_3$ is trivial). 
This configuration is depicted in Figure \ref{fig: good barbell quiver}, where we put $q_1 = \gamma_r\cdots\gamma_1$, $q_2 = \delta_s\cdots\delta_1$, $q_3 = \mu_t\cdots\mu_1$ and $q_4 = \nu_m\cdots\nu_1$. We consider the list
$$L:= (s(\gamma_2), \ldots, s(\gamma_r), s(\delta_1), \ldots, s(\delta_s), s(\beta), s(\mu_1), \ldots, s(\mu_t), s(\nu_1), \ldots, s(\nu_m))$$
of vertices. If there is repetition in $L$, take  $x$ to be such a repetition and assume that $x$ is minimal with the property that between the first two occurrences of $x$ in the ordered list $L$, there is no other repeated vertices. As discussed below, among the following cases, only case (3)(ii) results in the construction of a good barbell subquiver.

\begin{figure}
    \centering
\begin{tikzpicture}

 \draw [->] (0,0) --(1.7,0);
    \node at (1,0.2) {$\alpha$};
    \node[blue] at (1.75,0) {$\bullet $};
        \node[blue] at (1.75,0.25) {$a$};
 \draw [->] (1.75,0)--(1.25,-0.6);
 \draw [dashed] (1.2,-0.65)--(1.2,-1.2);
 \draw [->] (1.2,-1.25)--(1.2,-2);
    \node[red] at (1.2,-2.1) {$\bullet$};
        \node [red] at (0.95,-2) {$a'$};
     \draw [->] (1.15,-2.1)--(0.5,-2.7);
     \node at (1,-2.5) {$\gamma_1$};
\draw [dashed] (0.5,-2.7) to [bend right=100] (2,-2.7);
     \draw [->] (2,-2.7)--(1.26,-2.1);  
     \node at (1.5,-2.5) {$\gamma_r$};
     \draw [->] (1.26,-2.05)--(2.1,-1.5); 
     \node at (1.8,-1.9) {$\delta_1$};
     \draw [dashed] (2.1,-1.5) -- (2.1,-0.6);
       \draw [->] (2.1,-0.6)--(1.8,-0.05);
       \node at (2.1,-0.25) {$\delta_s$};
 \draw [->] (1.8,0) --(4.42,0);
 \node at (3,0.2) {$\beta$};
 \node[blue] at (4.5,0) {$\bullet$};
   \node[blue] at (4.5,0.25) {$b$};
 

 \draw [->] (4.5,0) --(4,-0.6);
 \node at (4.1,-0.25) {$\mu_1$};
  \draw [dashed] (4,-0.6) --(4,-1.25);
    \draw [->] (4,-1.25) --(4,-2);
    \node at (3.82,-1.5) {$\mu_t$};
    \node[red] at (4,-2.1) {$\bullet$};
    \node [red] at (3.75,-2) {$b'$};
    \draw [->] (3.95,-2.1)--(3.4,-2.7); 
    \node at (3.87,-2.5) {$\nu_1$};
\draw [dashed] (3.4,-2.7) to [bend right=100] (4.9,-2.7);
     \draw [->] (4.9,-2.7)--(4.05,-2.1);  
         \node at (4.42,-2.5) {$\nu_m$};
     \draw [->] (4.05,-2.05)--(4.9,-1.5); 
     \draw [dashed] (4.9,-1.5) -- (4.9,-0.6);
       \draw [->] (4.9,-0.6)--(4.55,-0.05);
\end{tikzpicture}
    \caption{Existence of a good barbell subquiver}
    \label{fig: good barbell quiver}
\end{figure}

\medskip

\emph{Case $(1)$}: We have $x=s(\gamma_i)$ for $i \ge 2$.  In this case, $p$ contains two parallel paths $p', p''$ starting at $a'$ and ending at vertex $x$ and not sharing any arrows or vertices (other than starting and ending points). Indeed, we take $p' = \gamma_{i-1} \cdots \gamma_{1}$ where $x = s(\gamma_i)$, and $p''$ start with $\delta_1$ and minimal with the properties of being good and ending at vertex $x$. Being subpaths of $p$, none of them belong to $I$. We observe that we are in the setting of Lemma \ref{LemmaParallel}, where $c = a', d = x$ and $p'$ and $ p''$ respectively play the role of $p_1$ and $ p_2$ of that lemma. By Lemma \ref{LemmaParallel}, there is an arrow of $p'$ or $ p''$ ending in $a'\ne x$. This contradicts the assumption that $p'$ or $ p''$ share no vertices other than their starting and ending points. 

\medskip

\emph{Case $(2)$}: We have $x=s(\nu_j)$ for $j \ge 2$ which is the second occurrence of $x$ in $L$. In this case, $p$ contains two parallel paths $p', p''$ ending at $b'$ and starting at vertex $x$ and not sharing any arrows or vertices (other than starting and ending points). Indeed, we take $p' = \nu_m \cdots \nu_j$ and $p''$ ending with $\mu_t$ (or $\beta$ if $q_3$ is trivial) and being minimal with the properties of being good and starting at $x$. This leads to a similar contradiction as in the previous case.

\medskip

\emph{Case $(3)$}: We are not in cases $(1)$ or $(2)$. Thus, the repetition $x$ involves vertices in the list $$L':=\{a'=s(\delta_1), \ldots, s(\delta_s), s(\beta), s(\mu_1), \ldots, s(\mu_t), s(\nu_1)=b'\}.$$ We recall that $x$ was chosen to be minimal. Hence, the vertex $x$ is such that between the two occurrences of $x$ in $L'$, we get no other repetitions. We separate this into some subcases.

\medskip

\emph{Subcase $(i)$}: Assume first that $x=a'$ is the repeated vertex. Note that $q_2$ does not revisit $a'$ by definition, so that the second time $x$ appears in $L'$ is $x = e(\mu_i)$ for some $i$ or $x=b$. Now, a good path ending at $a'$ either ends with $\gamma_r$ or with the last arrow (call it $\epsilon$) of $p_{11}$.  Note that if $p_{11}$ is trivial, then $a=a'$ and we take $\epsilon = \alpha$.  We also observe that $\gamma_1 \gamma_r$ is bad while $\gamma_1 \epsilon$ is good, therefore $\gamma_r \ne \epsilon$.

Consider first the case where $x = e(\mu_i)$. Then $\mu_i = \gamma_r$ or $\mu_i = \epsilon$. In the first case, we get a contradiction on the minimality of $x=a'$. To see this, if $\mu_i = \gamma_r$, then $\gamma_r$ cannot be a loop, because $q_3$ does not revisit vertices. Hence, in this case, $s(\gamma_r)$ belongs to $L$ and the repetition of it occurs in $L$ before the second occurrence of $x=a'$, which is a contradiction to the minimality of $x$ in $L$. 
Now, assume $\mu_i = \epsilon$. This implies that either $\mu_i\cdots\mu_1\beta$ is a subpath of $p_{11}$ or that $p_{11}$ is a subpath of $\mu_i\cdots\mu_1\beta$ (the shorter one is a subpath of the longer one). The first case cannot happen since $\beta$ does not appear in $p_1$ by definition. The second case leads to a contradiction of the minimality of $x=a'$, since the vertex $a$ repeats and the two occurrences appear before $e(\mu_i)$. 
Here, we note that $p_1$ does not contain $\beta$. Hence, if $p_{11}$ is a subpath of $\mu_i \dots \mu_1\beta$, then $p_{11}$ is a subpath of $\mu_i \dots \mu_1$, which implies that $a = s(\mu_j)$, for some $j$, so that $a = s(\beta) = s(\mu_j)$ appears between the two occurrences of $x=a'$ in the list $L$, a contradiction.

The remaining case is when $x=a'=b$. This means that $\beta \in \{\epsilon, \gamma_r\}$, which is impossible since $p_1\alpha$ does not contain $\beta$.

\medskip

\emph{Subcase $(ii)$}: Assume that $x \not \in \{a', b'\}$ is the repeated vertex. 
Then we may shorten $q_4q_3\beta q_2q_1$ and create a smaller good barbell quiver. In particular, as shown in Figure \ref{fig: Subcase (iii)}, we get a good barbell subquiver, where $C_L = \gamma_r \cdots \gamma_1$ and $C_R=\mu_j \cdots \mu_1 \beta \delta_s \cdots \delta_{i}$, and the bar is $\delta_{i-1} \cdots \delta_1$. Observe that in the new good barbell quiver, which is depicted with thick arrows and dashed segments, the arrow $\beta$ is part of a cycle. In particular, it is no longer in the bar.

\medskip

\emph{Subcase $(iii)$}: We have that $x=b'$ is the repeated vertex. We note that the first appearance of $x$ has to be at one of $s(\delta_2), \ldots, s(\delta_s), a$. Moreover, $s \ge 1$ as otherwise, we would be in subcase $(i)$ with $a=a' = b'$. Note also that no $\delta_i$ could be equal to $\mu_t$, as otherwise, the vertex $s(\delta_i)=s(\mu_t)$ would appear twice in $L'$, between the two occurrences of $x=b'$, which clearly contradicts the minimality of $x$. Also, no $\delta_i$ could be equal to $\beta$. Hence there is some $i$ such that $\delta_i = \nu_m$. This means that either $\delta_i \cdots \delta_1 \gamma_r \cdots \gamma_1$ is a subpath of $\nu_m \cdots \nu_1$ or the other way around. Since $\gamma_r \cdots\gamma_1$ and $\nu_m \cdots \nu_1$ are simple cycles and since there is no repeated vertices in $\delta_s \cdots \delta_1$, this means that $q_2$ is trivial and $\gamma_r \cdots \gamma_1 = \nu_m \cdots \nu_1$, implying that $x = a=a' = b'$, which is a case that was treated in Subcase (i).

\begin{figure}
    \centering
\begin{tikzpicture}

 \draw [->] (0,0) --(1.7,0);
    \node at (1,0.2) {$\alpha$};   
    \node[blue] at (1.75,0) {$\bullet $};
        \node[blue] at (1.75,0.25) {$a$};
 \draw [->] (1.75,0)--(1.25,-0.6);
 \draw [dashed] (1.2,-0.65)--(1.2,-1.2);
 \draw [->] (1.2,-1.25)--(1.2,-2);
    \node[red] at (1.2,-2.1) {$\bullet$};
        \node [red] at (0.95,-2) {$a'$};
     \draw [line width=0.5mm, ->] (1.15,-2.1)--(0.5,-2.7);
     \node at (1.03,-2.5) {$\gamma_1$};
\draw [line width=0.5mm, dashed] (0.5,-2.7) to [bend right=100] (2,-2.7);
     \draw [line width=0.5mm, ->] (2,-2.7)--(1.26,-2.1);  
     \node at (1.5,-2.5) {$\gamma_r$};

     \draw [line width=0.5mm, ->] (1.26,-2.05)--(1.8,-1.7); 
     \node at (1.7,-2.03) {$\delta_1$};

\draw [line width=0.5mm, dashed] (1.8,-1.7) -- (2.4,-1.3);
\draw [line width=0.5mm, ->] (2.4,-1.3)--(2.8,-1.05);
\node at (2.73,-1.45) {$\delta_{i-1}$};
       \node [red] at (2.9,-1) {$\bullet$};
       \node [red] at (2.9,-1.2) {$x$};
\draw [line width=0.5mm, ->] (2.9,-1)--(2.5,-0.7);
\node at (2.85,-0.62) {$\delta_{i}$};
\draw [line width=0.5mm, dashed] (2.5,-0.7) -- (2.1,-0.35);
\draw [line width=0.5mm, ->] (2.1,-0.35)--(1.8,-0.05);
\node at (2.2,-0.2) {$\delta_s$};
 \draw [line width=0.5mm, ->] (1.8,0) --(4.42,0);
 \node at (3,0.2) {$\beta$};
        \node[blue] at (4.5,0) {$\bullet$};
        \node[blue] at (4.5,0.25) {$b$};

 \draw [line width=0.5mm, ->] (4.5,0) --(4.2,-0.6);
 \node at (4.15,-0.25) {$\mu_1$};
  \draw [line width=0.5mm, dashed] (4.2,-0.57) --(3.5,-0.8);
   \draw [line width=0.5mm, ->] (3.5,-0.8) --(2.95,-1);
\node at (3.3,-0.64) {$\mu_{j}$};

  \draw [dashed] (2.9,-1) --(4,-1.2);

    \draw [->] (4,-1.25) --(4,-2);
    \node at (3.82,-1.5) {$\mu_t$};
    \node[red] at (4,-2.1) {$\bullet$};
    \node [red] at (3.75,-2) {$b'$};
    \draw [->] (3.95,-2.1)--(3.4,-2.7); 
    \node at (3.87,-2.5) {$\nu_1$};
\draw [dashed] (3.4,-2.7) to [bend right=100] (4.9,-2.7);
     \draw [->] (4.9,-2.7)--(4.05,-2.1);  
         \node at (4.42,-2.5) {$\nu_m$};
     \draw [->] (4.05,-2.05)--(4.9,-1.5); 
     \draw [dashed] (4.9,-1.5) -- (4.9,-0.6);
       \draw [->] (4.9,-0.6)--(4.55,-0.05);
\end{tikzpicture}
    \caption{Subcase (iii) of Case (3)}
    \label{fig: Subcase (iii)}
\end{figure}

\medskip

Having treated all cases, we therefore see that $q_4q_3\beta q_2q_1$ contains a subpath having $\beta$ which forms a good barbell subquiver. In particular, in the case where there is no repetition in $L$, then $q_4q_3\beta q_2q_1$ is of the required form. Otherwise, $(3)(ii)$ is the only possible case, where the construction is given.
\end{construction}

The construction above leads to the following proposition, which finishes our argument about the description of those non-hereditary biserial algebras which are minimal brick-infinite. 

\begin{proposition}
Let $A = kQ/I$ where $Q$ is biserial, a bisection is given and $I$ contains the generators as in Theorem \ref{DescriptionBiserial}. If A is minimal brick-infinite and $Q$ contains a good barbell subquiver, then $A$ itself is a barbell algebra. In particular, $A$ is gentle.
\end{proposition}

\begin{proof}
By Construction \ref{Remark: Existence of good barbell subquiver}, there exists a good barbell subquiver $Q'$ in $(Q,I)$, with a unique maximal good path, say $w$ of $Q'$, which could be seen as a good path of $Q$ and it is \emph{a priori} not maximal in $Q$. Assume $w$ in $Q$ starts at $a$ and ends at $b$.
By Lemma \ref{MonomialRelationsBiserial}, the ideal $I$ can be generated by the relations $R_{\beta\alpha}$ as in Theorem \ref{DescriptionBiserial} for the bad paths $\beta\alpha$ of length two, plus possibly some monomials relations, which can be taken to be good paths. Let $J$ be the ideal of $kQ$ generated by the following monomial relations.   
\begin{enumerate}
    \item The vertices and arrows not appearing in $Q'$;
    \item The monomial relations from $I$;
    \item For each bad path $\nu\mu$ of length two with corresponding relation $R_{\nu\mu} = \nu\mu - \lambda_{\nu\mu} r_{\nu\mu} \mu$, with $\lambda_{\nu\mu} \ne 0$, the relations $\nu\mu$ and $r_{\nu\mu} \mu$ whenever they are paths of $Q'$.
\end{enumerate}
We observe that $J$ contains $I$. Clearly, $kQ/J$ is a string algebra. Let $q$ be a good path in $Q'$, where the bisection is inherited from that of $Q$. 

We claim that $q$ is not in $J$. It follows from the definition of good barbell subquiver that $q$ does not contain any monomial relation from $I$. Hence, to prove the claim, it is sufficient to show that if $\nu\mu$ is a bad path in $(Q,I)$, then $q$ is not  equal to the good path $r_{\nu\mu}\mu$ in $(Q,I)$, where $r_{\nu\mu}$ is such that $R_{\nu\mu} = \nu\mu - \lambda_{\nu\mu} r_{\nu\mu} \mu$ for a scalar $\lambda_{\nu\mu}$. 
Assume otherwise, that is, suppose there exists a bad path $\nu\mu$ such that $q = r_{\nu\mu}\mu$ as above.
First, we observe that $q$ is an ending subpath of $w$ in $Q'$; otherwise, the assumption $q = r_{\nu \mu}\mu$ implies that $w = \lambda_{\nu \mu}^{-1}q'\nu \mu q''$ in $(Q, I)$ for some paths $q', q''$ with $q'$ non-trivial, which in turns implies that $w \in I$ by Proposition \ref{BadPathSocle}(1) – a contradiction. 
Now, we have that $q$ is an ending subpath of $w$, in particular, it ends at vertex $b$. 
But since $Q'$ already has two arrows ending at $b$ and since $q$ cannot end with $\nu$, by part (2) of Theorem \ref{DescriptionBiserial}, we conclude that $\nu$ is the arrow on the bar (of the good barbell subquiver obtained from $Q$) ending at $b$. On the other hand, it follows from the definition of a good barbell quiver, particularly because $\nu$ belongs to the bar, that there is only one arrow ending with $s(\nu)$ in $Q'$. Since $\mu$ lies in $Q'$, we conclude that the bad path $\nu\mu$ lies on the good path $w$, which is a contradiction. This proves our claim that $q$ is not in $J$.

We have now proved that $kQ/J$, which is a quotient of $A$, is isomorphic to the algebra $kQ' / (kQ' \cap J)$ of a generalized barbell subquiver with a non-trivial bar. Since $A$ is minimal brick-infinite and because $kQ/J$ is brick-infinite, we get $A = kQ/J$ is itself a generalized barbell algebra.
\end{proof}

The preceding proposition, along with our assumptions in this section, completes our proof of Theorem \ref{Thm: Classification of min-brick-infinite biserials}. Consequently, we get the following result which gives a full classification of biserial algebras with respect to brick-finiteness.

\begin{corollary}\label{Cor: full classification of brick-finiteness of biserial algs}
A biserial algebra is brick-infinite if and only if it has a gentle quotient algebra $A=kQ/I$ such that $(Q,I)$ admits a band.
\end{corollary}

We note that our results also apply to some other families of algebras which have been studied in the literature, such as the weighted surface algebras introduced in \cite{ES}, as well as the stably biserial algebras, studied in \cite{Po}.
In particular, in the assertion of the preceding corollary, one can replace biserial algebras with weighted surface algebras or stably biserial algebras. The argument is quite straightforward but requires some considerations, which we leave to the interested reader.

\section{Some applications and problems}\label{Section:Some applications and problems}
Here we consider some consequences of our results in the preceding sections and propose a new treatment of brick-infinite algebras. In doing so, we view some classical results through a new lens which better motivates some questions posed below. 
As before, we work over an algebraically closed field and, unless specified otherwise, $\Lambda$ denotes a tame algebra.
Also, recall that $G$ is a \emph{generic brick} if it is generic and $\End_{\Lambda}(G)$ is a division algebra. As in \cite{Ri2}, one can treat generic bricks as certain points of the spectrum of $\Lambda$. This is a conceptual generalization of spectrum of commutative rings to any arbitrary ring, first introduced by P. Cohn \cite{Co}. However, here we primarily study them from the algebraic and geometric viewpoints. For further details on spectrum of algebra, see \cite{Ri2}.

\subsection{Generic bricks and generic-brick-domestic algebras}
The Tame/Wild dichotomy theorem of Drozd \cite{Dr} plays a decisive role in the study of rep-infinite algebras. The family of tame algebras further refines into three disjoint subfamilies: domestic algebras, algebras of polynomial growth, and algebras of non-polynomial growth (for definitions and background, see \cite{Sk}).
To motivate a modern analogue of domestic algebras, we recall a fundamental theorem of Crawley-Boevey which gives a conceptual characterization of tameness, as well as domestic algebras.
Following \cite{CB1}, we say that $\Lambda$ is \emph{generically tame} if for each $d \in \mathbb{Z}_{>0}$, there are only finitely many (isomorphism classes) of generic modules of endolength $d$. 

\begin{theorem}[\cite{CB1}]\label{Thm Crawley-Boevey tame/domestic alg.}
An algebra $\Lambda$ is tame if and only if it is generically tame. Moreover, $\Lambda$ is domestic exactly when it admits only finitely many isomorphism classes of generic modules.
\end{theorem}

As we do henceforth, the characterization of tame and domestic algebras in the above theorem can be adopted as their definition. Moreover, for a $m \in \mathbb{Z}_{\geq 0}$, the algebra $\Lambda$ is $m$-domestic if and only if it admits exactly $m$ generic modules (see \cite[Corollary 5.7]{CB1}). 
Analogously, we say $\Lambda$ is \emph{$m$-generic-brick-domestic} if it admits exactly $m$ (isomorphism classes) of generic bricks. In general, we call $\Lambda$ \emph{generic-brick-domestic} if it admits only finitely many generic bricks.

It is known that $\Lambda$ is rep-finite if and only if $\Ind(A)\setminus \ind(A)=\emptyset$ (see \cite{Au}). Furthermore, by \cite{CB1} this is equivalent to the non-existence of generic module. Hence, the family of $0$-domestic algebras is the same as that of rep-finite algebras.
Moreover, Theorem \ref{Thm:Sentieri} implies that any brick-finite algebra is $0$-generic-brick-domestic. However, we do not know whether the converse is true in general.
Thanks to our new results, we can affirmatively answer this question for the family of biserial algebras and further conjecture that this holds in general.

\begin{remark}
We note that, in contrast to the notion of generic-brick-domestic algebras defined above, one can call a tame algebra \emph{$m$-brick-domestic} algebra if for any $d\in \mathbb{Z}_{>0}$, there are most $m$ one-parameter families of bricks of length $d$. We observe that this notion is different from $m$-generic-brick-domestic algebras and further studying of connections between these two notions could be interesting. We do not treat this comparison in this paper.
\end{remark}

Before restating our main conjecture for arbitrary algebras, we recall some basic notions and facts on the representation varieties of algebras. Recall from Section \ref{Section:Preliminaries} that $\Lambda$ is brick-discrete if for every $\mathcal{Z}$ in $\Irr(\Lambda)$, the set $\brick(\mathcal{Z})$ contains at most one brick (up to isomorphism).
If $\brick(\mathcal{Z}) \neq \emptyset$, then $\mathcal{Z}$ is called a \emph{brick component}. Because $\brick(\mathcal{Z})$ is always an open subset of $\mathcal{Z}$, each brick component is an indecomposable component.
Thus, $\Lambda$ is brick-discrete exactly when each brick component in $\Irr(\Lambda)$ is of the form $\mathcal{Z}=\overline{\mathcal{O}}_M$, for some $M$ in $\brick(\Lambda)$.
As in the Introduction, $\Lambda=kQ/I$ is called \emph{brick-continuous} if it is not brick-discrete. That is, there exists $\underline{d} \in \mathbb{Z}^{Q_0}_{\geq 0}$ and $\mathcal{Z}$ in $\Irr(\Lambda,\underline{d})$ such that $\brick(\mathcal{Z})$ contains infinitely many orbits of bricks.
In \cite{Mo2}, the first-named author conjectured that an algebra is brick-finite if and only if it is brick-discrete (see also Conjecture \ref{My Old Conjecture}).
Below, we propose a stronger version of this conjecture which also implies Theorem \ref{Thm:Sentieri}.

\begin{conjecture}\label{Conjecture for arbitrary algs}
For any algebra $\Lambda$, the following are equivalent:
\begin{enumerate}
    \item $\Lambda$ is brick-infinite;
    \item $\Lambda$ is brick-continuous;
    \item $\Lambda$ admits a generic brick.
\end{enumerate}
\end{conjecture}

In the preceding conjecture, observe that $(2)$ evidently implies $(1)$, and from Theorem \ref{Thm:Sentieri} it is immediate that $(3)$ implies $(1)$.
Furthermore, the implication $(3)$ to $(2)$ holds if there exist a generic brick which satisfies the assumption of the next proposition. In particular, this condition always holds for tame algebras.

\begin{proposition}[\cite{Ri2}]\label{Prop: Ringel on generic bricks}
Let $\Lambda$ be an algebra and $G$ be a generic brick such that $\End_{\Lambda}(G)$ is finitely generated over its center. Then, $G$ gives rise to a one-parameter family of bricks in $\brick(\Lambda)$.
\end{proposition}

To verify Conjecture \ref{Conjecture for arbitrary algs} for the family of tame algebras, it suffices to show the implication $(1)\Rightarrow(3)$ and for that one can reduce to minimal brick-infinite tame algebras.
In the following theorem, we consider certain families of such algebras and extend a classical result of Ringel on the hereditary case.
We remark that if $G_1$ and $G_2$ are two non-isomorphic generic modules in $\Ind(\Lambda)$, they induce two distinct $1$-parameter families of indecomposable modules in $\ind(\Lambda)$ (see \cite{CB1}). As mentioned earlier, the next theorem can be stated in the language of spectrum of algebras, as in \cite{Ri2}.

\begin{theorem}\label{Thm:Generarlization of Ringel's 1-generic-brick-domestic}
Let $\Lambda$ be a minimal brick-infinite tame algebra. If $\Lambda$ is hereditary or biserial, then
\begin{itemize}
    \item $\Brick(\Lambda)$ has a unique generic brick;
    \item $\brick(\Lambda)$ is the disjoint union of an infinite discrete family with $\{X_{\lambda}\}_{\lambda \in k*}$, where all $X_{\lambda}$ are of the same dimension.
\end{itemize}
In particular, in either of these cases $\Lambda$ is brick-continuous and $1$-generic-brick-domestic.
\end{theorem}

\begin{proof}
If $\Lambda$ is hereditary, it is the path algebra of some $\widetilde{A}_n, \widetilde{D}_m, \widetilde{E}_6, \widetilde{E}_7$ or $\widetilde{E}_8$, where $n\in \mathbb{Z}_{\geq 1}$ and $m\in \mathbb{Z}_{\geq 4}$. In this case, the assertions follow from the main result of \cite{Ri2}.
If we assume $\Lambda=kQ/I$ is biserial and non-hereditary, Theorem \ref{Thm: Classification of min-brick-infinite biserials} implies that $\Lambda$ is a generalized barbell algebra. As shown in \cite{Mo1}, every generalized barbell algebra admits a unique band $w$ for which the band module $M(w,\lambda)$ is a brick, for all $\lambda \in k^*$. In particular, explicit description of $w$ depends on the length and orientation of the bar $\mathfrak{b}$ in the generalized barbell quiver $(Q,I)$, as depicted in Figure \ref{fig:Generalized barbell quiver}. 
If $l(\mathfrak{b})$ denotes the length of $\mathfrak{b}$, we need to consider the two cases $l(\mathfrak{b})>0$ and $l(\mathfrak{b})=0$, as discussed below.
In the following, by $C_L$ and $C_R$ we denote respectively the left and right cyclic strings in $(Q,I)$ and assume $C_L=\alpha \nu^{\epsilon_p}_{p} \cdots \nu^{\epsilon_2}_2 \beta$ and $C_R=\gamma \mu^{\epsilon'_q}_q \cdots \mu^{\epsilon'_2}_2 \delta$, with $\mu_j, \nu_i \in Q_1$ and $\epsilon_i, \epsilon'_j \in \{\pm 1\}$, for every $1 \leq i \leq p$ and $1 \leq j \leq q$.

If $l(\mathfrak{b})>0$, let $s(\mathfrak{b})=x$ and $e(\mathfrak{b})=y$, and suppose $\mathfrak{b}=\theta^{\epsilon_d}_d \cdots \theta^{\epsilon_2}_2 \theta^{\epsilon_1}_1$ with $\theta_i \in Q_1$, for all $1 \leq i \leq d$. Without loss of generality, we can assume $\epsilon_1=1$, because the case $\epsilon_1=-1$ is similar. Then, by \cite[Proposition 5.6]{Mo1}, $w:= \mathfrak{b}^{-1} C_R \mathfrak{b} C_L$ gives us the desired band in $(Q,I)$, which we use to construct a generic brick $G=(\{G_i\}_{i\in Q_0}, \{G_{\eta}\}_{\eta \in Q_1})$ over $\Lambda$.
Starting from $x$, put a copy of $k(t)$ at $i\in Q_0$ each time $w$ passes through $i$. As the result, for each vertex $i$ that belongs to $\mathfrak{b}$ we have $G_i=k(t)\oplus k(t)$, whereas at the remaining vertices $j$ we get $G_j=k(t)$. As for the linear maps, as we go through $w$, for all arrows $\eta \in Q_1$ except for the second occurrence of $\theta_1$, we put the identity map between the two copies of $k(t)$ and in the direction of $\eta$, whereas the map from the first copy of $k(t)$ at $x$ to the second copy of $k(t)$ at $e(\theta_1)$ is given by multiplication by $t$.
Then, an argument similar to \cite[Proposition 5.6]{Mo1} shows that $\End_{\Lambda}(G)\simeq k(t)$, and from the construction it is clear $G$ is of finite dimension over $k(t)$. Hence, $G$ is the desired generic brick. 

Note that if $l(\mathfrak{b})=0$, the strings $C_L$ and $C_R$ cannot be serial simultaneously (otherwise $\Lambda$ will be infinite dimensional). In this case the desired band is given by $w:= C_R C_L$ and an argument similar to the above case gives the explicit construction of the generic brick. 

To show the uniqueness of this generic brick $G$, assume otherwise and let $G'$ be a generic brick in $\Brick(\Lambda)$ which is not isomorphic to $G$. By Proposition \ref{Prop: Ringel on generic bricks}, both $G$ and $G'$ induce $1$-parameter families of bricks in $\brick(\Lambda)$, say respectively $\{X_{\lambda}\}_{\lambda \in k^*}$ and $\{Y_{\lambda}\}_{\lambda \in k^*}$. From \cite{BR}, all these bricks are band modules.
Moreover, from \cite{CB1}, we know that these two $1$-parameter families are distinct, which implies there must come from two distinct bands in $(Q,I)$, say $w$ and $w'$, for which $M(w,\lambda)$ and $M(w',\lambda)$ are bricks. This contradicts the uniqueness of $w$, as shown in \cite[Proposition 5.6 and Proposition 7.10]{Mo1}.

Finally, observe that each $X$ in $\brick(\Lambda)$ is either a band module of the form $M(w,\lambda)$, for some $\lambda \in k^*$, or else is a string module. The former type gives a $1$-parameter family, whereas string bricks form a countable (discrete) infinite family. In particular, if $l(\mathfrak{b})>0$, we note that each string module $M(w^d)$ is a brick, where $d \in \mathbb{Z}_{\geq 1}$ (for details, see \cite{Mo1}). If $l(\mathfrak{b})=0$, the explicit description of an infinite family of string modules which are bricks is given in \cite[Proposition 7.10]{Mo1}.
\end{proof}

\begin{figure}
    \centering
\begin{tikzpicture}

 \draw [->] (1.25,0.75) --(2,0.1);
    \node at (1.7,0.55) {$\alpha$};
 \draw [<-] (1.25,-0.75) --(2,0);
    \node at (1.7,-0.5) {$\beta$};
  \draw [dashed] (1.25,0.75) to [bend right=100] (1.25,-0.75);
   \node at (1.3,0) {$C_L$};
    \node at (2,0) {$\bullet $};
    \node at (2.1,-0.2) {$x$};
 \draw [dashed] (2,0) --(2.75,0);
 \node at (2.75,0) {$\bullet$};
 \draw [dashed] (2.75,0) --(4,0);
 \node at (4,0) {$\bullet$};
 \draw [dashed] (4,0) --(4.75,0);
 \node at (3.5,0.3) {$\mathfrak{b}$};
 
 \node at (4.75,0) {$\bullet$};
 \draw [<-] (5.5,0.75) --(4.75,0);
 \node at (5,0.45) {$\delta$};
 \draw [->] (5.50,-0.8) --(4.8,-0.05);
    \node at (5,-0.55) {$\gamma$};
  \draw [dashed] (5.55,0.8) to [bend left=100] (5.55,-0.8);
   \node at (5.5,0) {$C_R$};
   \node at (4.65,-0.2) {$y$};
\end{tikzpicture}

    \caption{Generalized barbell quiver}
    \label{fig:Generalized barbell quiver}
\end{figure}

We remark that (generalized) barbell quivers are of non-polynomial growth (see \cite{Ri1}). That means, roughly speaking, as long as the behavior of all indecomposable modules is concerned, generalized barbell algebras are among the most complicated type of tame algebras. However, the preceding theorem implies they are always \emph{$1$-generic-brick-domestic}. Hence, with respect to this modern criterion, generalized barbell algebras are among the simplest type of brick-infinite algebras (see also Subsection \ref{Subsection: Domestic vs. generic-brick-domestic}).

As a consequence of the preceding theorem, we get the following result. In particular, this proves Corollary \ref{Corollary-Introduction: all equivalences for brick-inf biserial algs}.

\begin{corollary}\label{Corollary: all equivalences for brick-inf biserial algs}
Let $\Lambda=kQ/I$ be a biserial algebra. The following are equivalent:

\begin{enumerate}
    \item $\Lambda$ is brick-infinite;
    \item $\Lambda$ is brick-continuous;
    \item $\Lambda$ admits a generic brick;
    \item There is an infinite family of non-isomorphic bricks of length $d \leq 2|Q_0|$;
    \item In a brick component $\mathcal{Z}$ in $\Irr(\Lambda)$, there a rational curve $\mathcal{C}$ of non-isomorphic bricks $\{M_{\lambda}\}$ such that $\mathcal{Z}=\overline{\bigcup_{\lambda\in \mathcal{C}}\mathcal{O}_{M_{\lambda}}}$;
    \item For some $\theta \in K_0(\Lambda)$, there exist infinitely many non-isomorphic $\Lambda$-modules which are $\theta$-stable.
\end{enumerate}
\end{corollary}

Before we present a proof, note that this corollary gives novel algebro-geometric realizations of $\tau$-tilting (in)finiteness of biserial algebras (see Theorem \ref{Thm:tau-tilting finiteness equivalences}). Viewed from this perspective, these results also extend some earlier work on the family of special biserial algebras (see \cite{STV}). 

\begin{proof}
First, we note that a tame algebra $\Lambda$ is brick-continuous if and only if for some $\mathcal{Z} \in \Irr(\Lambda)$ we have $\mathcal{Z}=\overline{\bigcup_{\lambda\in \mathcal{C}}\mathcal{O}_{X_{\lambda}}}$, where $\{X_\lambda\}_{\lambda \in \mathcal{C}}$ is a rational curve of bricks in $\mathcal{Z}$ (for details, see \cite{CC2}). This implies the equivalence of $(2) \Longleftrightarrow (5)$.
Furthermore, observe that the family of biserial algebras is quotient-closed, meaning that any quotient of a biserial algebra is again biserial. Hence, without loss of generality, we assume $\Lambda$ is a minimal brick-infinite biserial algebra. Then, by Theorem \ref{Thm:Generarlization of Ringel's 1-generic-brick-domestic}, $(1)$ implies $(3)$, from which we conclude the equivalences $(1)\Longleftrightarrow(2) \Longleftrightarrow(3)$ (see also the paragraph following Conjecture \ref{Conjecture for arbitrary algs}). 
Moreover, by the proof of Theorem \ref{Thm:Generarlization of Ringel's 1-generic-brick-domestic}, the unique generic brick $G$ on $\Lambda$ is of endolength $d \leq 2|Q_0|$, which induces an infinite family of band modules of length $d$ in $\brick(\Lambda)$. This shows $(3) \Rightarrow (4)$ and the reverse implication is immediate from $(4)\Rightarrow (1)$. Hence, the first five parts are equivalent.

It is well-known that if $M \in \modu(\Lambda)$ is $\theta$-stable, for some $\theta \in K_0(\Lambda)$, then $M$ is a brick. Hence, $(6)\Rightarrow (1)$ is immediate. To finish the proof, we note that $\Lambda$ admits an infinite subfamily of $\brick(\Lambda)$ consists of band modules (of the same length). These are known to be homogeneous and by a result of Domokos \cite{Do}, the are $\theta$-stable, for some $\theta \in K_0(\Lambda)$ (for explicit description of $\theta$, also see \cite[Lemma 2.5]{CKW}).
\end{proof}

\subsection{Domestic vs. generic-brick-domestic}\label{Subsection: Domestic vs. generic-brick-domestic}
In this subsection we highlight some fundamental differences between the two notions of domesticness for tame algebras. In particular, we present several examples to better clarify some important points and motivate some questions which could be further pursued. In doing so, we give specific attention to string algebras, because they provide a more tractable setting.

As remarked earlier, generalized barbell algebras give an explicit family of tame algebras which are not domestic but always $1$-generic-brick-domestic. This naturally raises the question if there are examples of tame algebras which are $n$-domestic and $m$-generic-brick-domestic, for arbitrary $m$ and $n$ in $\mathbb{Z}_{\geq 0}$. Evidently, we need to additionally assume $m\geq n$. Also, it is natural to ask whether there exist tame algebras which are $n$-generic-brick-domestic but not of polynomial growth.
In the following examples, we give an explicit algorithm for constructing such algebras and prove Corollary \ref{Cor: m-domestic and n-gen-brick-domestic}.

Before we give a set of examples, let us fix some notation and terminology that will be handy. We encourage the reader to observe the bound quiver given in Example \ref{Example: m-domestic algebras which are 0-generic-brick-domestic} while going through the following definition. Below, by $\Vec{A}_d$ we denote a linearly oriented copy of quiver $A_d$.

\begin{definition}\label{Def: A-nody gluing}
Let $\Lambda=kQ/I$ and $\Lambda'=kQ'/I'$ be two algebras and suppose $R$ and $R'$ are two minimal sets of generators, respectively, for the ideals $I$ and $I'$. 
Let $u$ be a sink in $Q$ and $w$ be a source in $Q'$. Furthermore, assume $\{\alpha_1, \alpha_2, \dots, \alpha_p \}$ be the set of all arrows in $Q$ ending at $u$, and $\{\gamma_1, \gamma_2, \dots \gamma_q\}$ be the set of all arrows in $Q'$ starting at $w$.

For any $d \in \mathbb{Z}_{> 0}$, the \emph{$\Vec{A}_d$-nody gluing} of $(Q,I)$ and $(Q',I')$ from $u$ to $w$ is the new quiver $(\widetilde{Q}, \widetilde{I})$ obtained as follows:

\begin{itemize}
    \item $\widetilde{Q}$ is the result of connecting $Q$ to $Q'$ via $\Vec{A}_d=\beta_{d-1}\dots\beta_1$ which begins at $u$ and ends at $w$.
    
    \item A generating set for the ideal $\widetilde{I}$ is given by 
    
    $R \cupdot R' \cupdot \{\beta_{s+1}\beta_s| 1\leq s\leq d-2\}\cupdot \{\beta_1\alpha_i| 1\leq i \leq p \} \cupdot \{\gamma_j\beta_{d-1}|1\leq j \leq q\}$.

\end{itemize}

In particular, all vertices of $\Vec{A}_d$, including those identified with $u$ and $w$, are nodes in $(\widetilde{Q}, \widetilde{I})$.
\end{definition}

In the previous definition, if $d=1$, then $(\widetilde{Q}, \widetilde{I})$ is obtained by identifying the vertices $u$ and $w$ and imposing all the monomial quadratic relations at every possible compositions $\gamma_j\alpha_i$, for all $1\leq i \leq p$ and $1\leq j \leq q$. Hence, the new vertex is evidently a node in the quiver $(\widetilde{Q}, \widetilde{I})$.

The following lemma is handy in discussing the following examples.

\begin{lemma}\label{Lem: generic modules over glued algebras}
With the same notion and convention as in Definition \ref{Def: A-nody gluing}, let $(\widetilde{Q}, \widetilde{I})$ be the $\Vec{A}_d$-nody gluing of two string algebras $\Lambda=kQ/I$ and $\Lambda'=kQ'/I'$.
Provided $d \in \mathbb{Z}_{\geq 2}$,
no generic module over $\widetilde{\Lambda}=k\widetilde{Q}/\widetilde{I}$ is supported on the arrows of $\Vec{A}_d$.
That is, every generic module in $\Ind (\widetilde{\Lambda})$ belongs to either $\Ind(\Lambda)$ or $\Ind(\Lambda')$.
\end{lemma}

Before we prove the lemma, let us remark that the assumption $\mathbb{Z}_{\geq 2}$ is essential (see Example \ref{Example: gentle algebra not generic-brick-domestic}).

\begin{proof}
From the construction of $(\widetilde{Q}, \widetilde{I})$, it is immediate that $\widetilde{\Lambda}$ is also a string algebra. Furthermore, all vertices of $\widetilde{Q}$ that belong to $\Vec{A}_d$ are nodes in $(\widetilde{Q}, \widetilde{I})$. That is, every composition of arrows in $\widetilde{Q}$ at these vertices falls in $\widetilde{I}$. Therefore, if $\widetilde{G}$ is a generic module over $\widetilde{\Lambda}$, it induces a $1$-parameter family of band modules. In particular, if $\widetilde{G}$ is supported on any of the arrows of $\Vec{A}_d$, so are the induced band modules. This is impossible, because there is no band in $(\widetilde{Q}, \widetilde{I})$ which is supported on an arrow of $\Vec{A}_d$. This completes the proof.
\end{proof}

In \cite{Ri1}, the author introduced a particular type of string algebras, so-called \emph{wind wheel} algebras, which played a crucial role in classification of minimal representation-infinite (special) biserial algebras. Since the definition of wind wheel algebras is technical, here we only recall the general configuration of their bound quivers and for their explicit constructions, we refer the reader to \cite{Ri1}.
In general, each wind wheel algebra has a bound quiver of the following form:

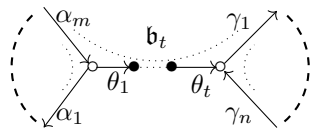
\begin{center}
\begin{tikzpicture}[scale=1]

\draw [->] (1,-0.8) --(0.3,-0.05);
    \node at (0.5,-0.65) {$\gamma_{n}$};
\draw [<-] (1,0.75) --(0.3,0.05);
     \node at (0.5,0.6) {$\gamma_{1}$};
    \draw [dotted] (0.7,-0.3) to [bend left=50](0.7,0.3);
    
    \draw [dashed, thick] (1.1,-0.72) to [bend right=50](1.1,0.72);

     \node at (0.25,0) {$\circ$};
       \node at (0.,-0.25) {$\theta_{t}$};
\draw [<-] (0.2,0) --(-0.35,0); 
    \node at (-0.4,0) {$\bullet$};
    \node at (-0.9,0) {$\bullet$};
\draw [dotted] (-0.9,0) --(-0.4,0); 
\draw [->] (-1.4,0) --(-0.95,0); 
    \node at (-1.45,0) {$\circ$};
    \node at (-1.1,-0.2) {$\theta_1$};
    \node at (-0.6,0.4) {${\mathfrak{b}_{t}}$};

\draw [->] (-2.1,0.8) --(-1.5,0.05);
 \node at (-1.7,0.65) {$\alpha_{m}$};

\draw [<-] (-2.1,-0.8) --(-1.5,0);
    \node at (-1.75,-0.65) {$\alpha_1$};
    \draw [dotted] (-1.85,-0.3) to [bend right=50](-1.85,0.3);
    
    \draw [dashed, thick] (-2.2,-0.72) to [bend left=50](-2.2,0.74);

    \draw [dotted] (-1.7,0.5) to [bend right=40](0.5,0.5);
    \end{tikzpicture}
\end{center}
where the bar ${\mathfrak{b}_{t}}$ is serial and of length $t \in \mathbb{Z}_{>0}$. Moreover, in addition to the monomial quadratic relations $\alpha_1\alpha_m$ and $\gamma_1\gamma_n$, we only have the monomial relation $\gamma_1 \theta_t\dots \theta_1\alpha_1$.
Note that the internal orientation of the left and right cycles are arbitrary, and this freedom is depicted by dashed segments. In particular, observe that we may have $\alpha_1=\alpha_m$ and similarly $\gamma_1=\gamma_n$.

\begin{example}[$m$-domestic algebras which are $0$-generic-brick-domestic]\label{Example: m-domestic algebras which are 0-generic-brick-domestic}
The wind wheel algebras are known to be $1$-domestic (see \cite{Ri1}). Moreover, in \cite{Mo1}, it is shown that they are always brick-finite. Hence, wind wheel algebras are $1$-domestic but $0$-generic-brick-domestic string algebras. 
Let $\Lambda=kQ/I$ and $\Lambda'=kQ'/I'$ be any pair of wind wheel algebras such that $Q$ has a sink and $Q'$ has a source. By description of their bound quivers, such sinks and sources cannot belong to their bars, thus they must be on one of their cycles that connect to vertices of degree $3$.
 
Then, consider an $\Vec{A}_3$-nody gluing of $(Q,I)$ and $(Q',I')$, such as the explicit example illustrated below. Note that, to distinguish the arrows of $\Vec{A}_3$ from those that belong to $Q$ and $Q'$, they are depicted by thicker solid arrows. Moreover, the solid arrows appearing on the left end (respectively, right end) of the large quiver $\widetilde{Q}$ belong to $Q$ (respectively, $Q'$). Moreover, those monomial quadratic relations which do not belong to $I$ nor $I'$ are imposed through the $\Vec{A}_3$-nody gluing, and they are shown with thicker dotted quadratic relation.

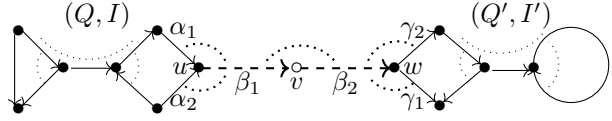
\begin{center}
\begin{tikzpicture}

\node at (3,3.2) {$(Q,I)$};

 \node at (1.95,3) {$\bullet$};
 \draw [->] (2,3) --(2.5,2.5); \node at (2.55,2.5) {$\bullet$};
 \draw [->] (2.5,2.5) --(2,2); \node at (1.95,1.95) {$\bullet$};
     \draw [<-] (1.9,2) --(1.9,2.95);
            \draw [dotted] (2.3,2.3) to [bend left=50](2.3,2.75);
  
\draw [->] (2.65,2.5) --(3.2,2.5); \node at (3.25,2.5) {$\bullet$};
            \draw [dotted] (2.25,2.9) to [bend right=50](3.5,2.9);

\draw [->] (3.25,2.5) --(3.75,3); \node at (3.8,3) {$\bullet$};
\draw [->] (3.85,2.97) --(4.3,2.55); \node at (4.35,2.5) {$\bullet$};
\node at (4.15,3) {$\alpha_1$};
\draw [<-] (4.35,2.5) --(3.85,2); \node at (3.8,1.95) {$\bullet$};
\node at (4.1,2.5) {$u$};
\node at (4.15,2.05) {$\alpha_2$};
\draw [->] (3.8,1.95) --(3.25,2.45);
            \draw [dotted] (3.45,2.3) to [bend right=50](3.45,2.75);
 \draw [dotted, line width=0.5mm] (4.7,2.6) to [bend right=50](4.15,2.77);
 \draw [dotted, line width=0.5mm] (4.7,2.4) to [bend left=50](4.15,2.27);

\draw [->, line width=0.5mm] (4.4,2.5) --(5.57,2.5); \node at (5.65,2.5) {$\circ$};

\node at (5,2.3) {$\beta_1$};
\node at (5.65,2.3) {$v$};
\node at (6.3,2.3) {$\beta_2$};

\draw [->, line width=0.5mm] (5.7,2.5) --(6.9,2.5); \node at (6.95,2.5) {$\bullet$};
             \draw [dotted, line width=0.5mm] (5.25,2.6) to [bend left=50](6.15,2.6);
 \draw [dotted, line width=0.5mm] (7.2,2.8) to [bend right=50](6.6,2.6);
 \draw [dotted, line width=0.5mm] (7.2,2.3) to [bend left=50](6.6,2.4);             
             
\node at (8.5,3.2) {$(Q',I')$};

\node at (7.2,2.5) {$w$};
\draw [->] (7,2.55) --(7.5,3); \node at (7.55,3) {$\bullet$};
\node at (7.2,2.05) {$\gamma_1$};
\draw [->] (7.6,2.97) --(8.1,2.55); \node at (8.15,2.5) {$\bullet$};
\node at (7.2,3) {$\gamma_2$};
\draw [->] (8.1,2.47) --(7.6,2.05); \node at (7.55,2) {$\bullet$};
\draw [<-] (7.5,2.05) --(7.02,2.47);
        \draw [dotted] (7.9,2.4) to [bend left=50](7.9,2.75);

\draw [->] (8.25,2.47) --(8.75,2.47); \node at (8.8,2.5) {$\bullet$};
            \draw [dotted] (7.85,2.9) to [bend right=50](8.95,3);

\draw[<-] (8.8,2.5) arc (185:545:0.5);
            \draw [dotted] (9,2.25) to [bend right=50](9,2.85);

\end{tikzpicture}
\end{center}

Let the above bound quiver $(\widetilde{Q}, \widetilde{I})$ be the $\Vec{A}_3$-nody gluing of $(Q,I)$ and $(Q',I')$. 
Then, $\widetilde{\Lambda}:=k\widetilde{Q}/\widetilde{I}$ is evidently a string algebra and it is easy to check that it is a $2$-domestic and brick-finite.
In particular, by Lemma \ref{Lem: generic modules over glued algebras}, every generic module in $\Ind(\widetilde{\Lambda})$ is a generic module over exactly one of the two wind wheel algebras, whose bound quivers are specified by solid arrows on the left and right ends of the above bound quiver.
Moreover, Corollary \ref{Cor: full classification of brick-finiteness of biserial algs} implies that $\widetilde{\Lambda}$ is brick-finite. Hence, $\widetilde{\Lambda}$ is $2$-domestic but $0$-generic-brick-domestic.

Through a recursive argument as above, for every $m \in \mathbb{Z}_{>0}$, one can explicitly construct $m$-domestic string algebras which are $0$-generic-brick-domestic.
\end{example}

In the next example, we consider some non-domestic tame algebras. In particular, we look at explicit string algebras that are of non-polynomial growth.

\begin{example}[Non-domestic tame algebras which are $0$-generic-brick-domestic]\label{0-generic-brick-domestic algebras which are not domestic} One can easily show that every (generalized) barbell algebra is of non-polynomial growth. In particular, it is a string algebra such that each arrow can appear in infinitely many distinct bands (for details, see \cite{Ri1}). Consider a barbell algebra, say $\Lambda=kQ/I$, such that $Q$ has at least one source, say $a$, and at least one sink, say $z$. It is well-known that by gluing these two vertices together, we obtain a new minimal representation-infinite algebra, say $\Lambda'=kQ'/I'$. That is, the vertices $a$ and $z$ are identified and all quadratic monomial relations are imposed at the new vertex (see \cite{Ri1}).
In particular, observe that $\Lambda'=kQ'/I'$ is again a string algebra and of non-polynomial growth.  However, the new vertex is a node and the new algebra $\Lambda'$ obtained via gluing vertices $a$ and $z$ becomes brick-finite (see \cite{Mo2}). In particular, it is $0$-generic-brick-domestic.
\end{example}

Below, we present another useful observation on the behavior of generic-brick-domestic algebras. In particular, the next example shows that there are tame algebras which are not generic-brick-domestic.

\begin{example}[Gentle but not generic-brick-domestic algebras]\label{Example: gentle algebra not generic-brick-domestic}
Consider the algebra $\Lambda=kQ/I$, where $(Q,I)$ is the following bound quiver, where all relations are quadratic. Note that $\Lambda$ is a gentle algebra, hence it is tame. We claim that $\Lambda$ is not generic-brick-domestic.

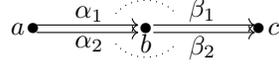
\begin{center}
\begin{tikzpicture}

 \node at (0.5,3) {$\bullet$};
  \node at (0.3,3) {$a$};
\draw [->] (0.55,3.05) --(1.9,3.05);
\draw [->] (0.55,2.95) --(1.9,2.95);
\node at (2,3) {$\bullet$};
\node at (2,2.8) {$b$};
 
  \node at (1.25,3.2) {$\alpha_1$};
  \node at (1.25,2.8) {$\alpha_2$};
 
   \draw [dotted] (1.6,2.8) to [bend right=50](2.4,2.8);
    \draw [dotted] (1.6,3.2) to [bend left=50](2.4,3.2);


 \node at (3.5,3) {$\bullet$};
  \node at (3.7,3) {$c$};
\draw [->] (2.1,3.05)--(3.45,3.05);
\draw [->] (2.1,2.95)--(3.45,2.95);

  \node at (2.75,3.25) {$\beta_1$};
    \node at (2.75,2.75) {$\beta_2$};

\end{tikzpicture}
\end{center}

First observe that $\Lambda$ admits two distinct quotient algebras which are hereditary, given by $\Lambda_a:=\Lambda/\langle e_a\rangle$ and $\Lambda_c:=\Lambda/\langle e_c\rangle$. Each of these quotient algebras is a path of the Kronecker quiver and admits a unique generic brick (see Theorem \ref{Thm:Generarlization of Ringel's 1-generic-brick-domestic}). 
Following the explicit description of morphisms between band modules, as in \cite{Kr}, one can show that the the band $w=\alpha_2^{-1}\beta_1^{-1}\beta_2\alpha_1$ is a brick band, hence gives rise to a generic brick.
More generally, for each $d \in \mathbb{Z}_{>0}$, consider the band $w_d:=\alpha_2^{-1}(\beta_1^{-1}\beta_2)^d\alpha_1$ in $(Q,I)$, where $(\beta_1^{-1}\beta_2)^d$ means the repetition of string $\beta_1^{-1}\beta_2$ for $d$ times. Evidently, if $d_1$ and $d_2 \in \mathbb{Z}_{>0}$ and $d_1 \neq d_2$, the bands $w_{d_1}$ and $w_{d_2}$ are distinct and the corresponding generic modules are non-isomorphic.
Again, using the same argument as above, one can show that each $w_d$ is a brick band, and therefore the corresponding generic module is a generic brick (see \cite{Kr}). Thus, the above gentle algebra is not generic-brick-domestic.
\end{example}

We finish this subsection with a proof of Corollary \ref{Cor: m-domestic and n-gen-brick-domestic}.
In particular, thanks to Theorem \ref{Thm:Generarlization of Ringel's 1-generic-brick-domestic} and the explicit constructions given in the preceding examples, we have the following result. For a string algebra $A$, we recall that if $A$ is not domestic, then $A$ must be of non-polynomial growth (see, for instance, \cite[Proposition 14.1]{Ri1} and its proof). Therefore, by \cite{CB1}, this particularly implies that $A$ admits infinitely many ismomorphism classes of generic modules.  Also, recall that a tame algebra is said to be generic-brick-domestic if it admits finitely many isomorphism classes of generic bricks. The following result asserts that even in the family of string algebras, a generic-brick-domestic algebra can be arbitrarily far from being domestic.

\begin{corollary}
Let $m$ and $n$ be a pair of non-negative integers. For each of the following mutually exclusive cases, there exists a string algebra $A$ such that
\begin{enumerate}
    \item $A$ is $m$-domestic and $n$-generic-brick-domestic, with $m\geq n$.
    
    \item $A$ is $n$-generic-brick-domestic but of non-polynomial growth.

    \item $A$ is not generic-brick-domestic, thus of non-polynomial growth.
\end{enumerate}
\end{corollary}
\begin{proof}
For case $(1)$, first observe that each rep-finite string algebra is $0$-domestic, therefore it is $0$-generic-brick-domestic. This treats the case $m=n=0$. For the case, $n=0$ but $m\in \mathbb{Z}_{>0}$, we refer to Example \ref{Example: m-domestic algebras which are 0-generic-brick-domestic}. 
We note that, if $\Lambda=kQ/I$ is a string algebra which is $m$-domestic and $n$-generic-brick-domestic, by an $\Vec{A}_2$-nody gluing of a wind wheel bound quiver to $(Q,I)$ we obtain a string algebra which is $(m+1)$-domestic but remains $n$-generic-brick-domestic (see Lemma \ref{Lem: generic modules over glued algebras} and Example \ref{Example: m-domestic algebras which are 0-generic-brick-domestic}). Moreover, via an $\Vec{A}_2$-nody gluing of a copy of $\widetilde{A}_d$ to the bound quiver to $(Q,I)$, we obtain a string algebra which is $(m+1)$-domestic and $(n+1)$-generic-brick-domestic.
Consequently, for each pair of integers $m$ and $n$ in $\mathbb{Z}_{>0}$, with $m\geq n$, we can construct a string algebra which is $m$-domestic and $n$-generic-brick-domestic.
This addresses all permissible choices of $m$ and $n$ in $(1)$.

Regarding case $(2)$, we first observe that each generalized barbell algebra is of non-polynomial growth and $1$-generic-brick-domestic (see Theorem \ref{Thm:Generarlization of Ringel's 1-generic-brick-domestic}). Similar to the construction in the Example \ref{Example: m-domestic algebras which are 0-generic-brick-domestic}, we can start from two generalized barbell algebras and via an $\Vec{A}_2$-nody gluing of the two copies we obtain a larger string algebra which is of non-polynomial growth and, by Lemma \ref{Lem: generic modules over glued algebras}, the new algebra is $2$-generic-brick-domestic. Alternatively, we can start from a generalized barbell algebra (with a source/sink in $C_L$ or $C_R$) and via an $\Vec{A}_2$-nody gluing with a copy of the Kronecker, we obtain a string algebra of non-polynomial growth which is $2$-generic-brick-domestic. 
By iteration, we can construct string algebras which are of non-polynomial growth, but $n$-generic-brick-domestic, for an arbitrary $n \in \mathbb{Z}_{\geq 0}$.

Finally, for case $(3)$, we refer to Example \ref{Example: gentle algebra not generic-brick-domestic}, where we construct an explicit example for a gentle algebra which is not generic-brick-domestic, and hence of non-polynomial growth.
\end{proof}

\subsection{Problems and remarks}
The rest of this section consists of some remarks and questions which outline our future work and new directions of research related to the scope of this paper.

\subsubsection*{Numerical condition for brick-finiteness}
For any algebra $\Lambda=kQ/I$ of dimension $d$, Bongartz has found some explicit numerical conditions for rep-finiteness and infiniteness of $\Lambda$ with respect to $d$ and some related constants (for details, see \cite[Section 7]{Bo2}). In particular, from \cite[Theorem 30]{Bo2} it follows that $\Lambda$ is rep-infinite if and only if for some $e \leq \max \{30,4d\}$, there is a an infinite family of (non-isomorphic) indecomposable $\Lambda$-modules of dimension $e$. This, by \cite{Sm}, then results in the existence of infinitely many $1$-parameter families of indecomposable modules asserted by the Second Brauer-Thrall conjecture (now theorem). 
This observation, along with our result on biserial algebras, naturally yields the question whether there is a similar criterion for brick-infiniteness of biserial algebras. Observe that, if $\Lambda=kQ/I$ is biserial, by Theorem \ref{Thm:Introduction:brick-infinite biserial algs}, $\Lambda$ is brick-infinite if and only if for some $t \leq 2|Q_0|$, there is a an infinite family of (non-isomorphic) bricks of dimension $t$. This, in particular, gives a very small and explicit bound which could be effectively applied to the study of $\tau$-tilting theory and stability conditions of biserial algebras, as in Corollary \ref{Corollary: all equivalences for brick-inf biserial algs}. 
More generally, one can ask the following question.

\begin{question}\label{Question: Numerical test for brick-finiteness}
Let $\Lambda$ be a $d$-dimensional algebra whose rank of Grothendieck group is $m$. Is there an explicit bound $b(d,m)$ in terms of $d$ and $m$ such that, to verify the brick-finiteness of $\Lambda$, it is sufficient to check if there is no infinite family of bricks of length less than $b(d,m)$?
\end{question}

Observe that, by definition, every domestic algebra is generic-brick-domestic, but Example \ref{0-generic-brick-domestic algebras which are not domestic} shows that the latter family is strictly larger.
Meanwhile, as shown in Example \ref{Example: gentle algebra not generic-brick-domestic}, there exist tame algebras which are not generic-brick-domestic.
Recently, Bodnarchuk and Drozd \cite{BD} introduced the brick-analogue of the classical notion of tameness, and called them brick-tame algebras. Furthermore, they give a new dichotomy theorem of algebras with respect to this criterion (for further details, see \cite{BD}). This notion is further studied by Carroll and Chindris \cite{CC1}. Note that strictly wild algebras are never brick-tame. It is a folklore open conjecture, attributed to Ringel, that the converse holds in general. Adopting this perspective in the study of bricks, combined with the characterization of tame/wild dichotomy through generic modules in \cite{CB1}, naturally raises the question whether there exists an analogue characterization of brick-tame/brick-wild algebras via generic bricks. 
The notion of brick-tameness allows one to approach brick-infinite algebras from a new perspective, particularly in the treatment of non-domestic tame algebras and those wild algebras which are not strictly wild.

As long as non-domestic tame algebras are concerned, a good knowledge of $m$-generic-brick-domestic algebras, for different $m \in \mathbb{Z}_{\geq 0}$, should provide fresh impetus to some classical problems.
In that sense, the following question could be of interest and relates to Theorem \ref{Thm:Generarlization of Ringel's 1-generic-brick-domestic}.

\begin{question}\label{Question:generic-brick-domesticness of tame min-brick-infinite algs}
Is every tame minimal brick-infinite algebra generic-brick-domestic?
\end{question}

As a question of the same nature, one can ask whether arbitrary minimal brick-infinite algebras always admit only finitely many generic bricks.

Finally, we note that for gentle algebras, recently Gei\ss{}, Labardini-Fragoso and Schr\"oer\cite{G+} have extensively studied the string and band components. Such results directly apply to the family of minimal brick-infinite biserial algebras, because they are always gentle (see Theorem \ref{Thm: Introduction: min-brick-inf biserial}). Hence, it is natural to further investigate that approach and try to extend the setting of \cite{G+}.

\vskip 0.6cm

\textbf{Acknowledgements.} The authors would like to thank William Crawley-Boevey for some helpful discussions related to the scope of this work. We also thank Aaron Chan for some interesting conversations about gentle algebras coming from surfaces.
Finally, we thank an unknown referee for several good comments that improved the exposition of this paper.

\end{document}